\documentclass[12pt,reqno]{amsart}
\usepackage{amsmath,latexsym,amssymb,verbatim,euscript}
\numberwithin{equation}{section}

\theoremstyle{definition}
\newtheorem{thm}{Theorem}
\newtheorem*{thm*}{Theorem}
\newtheorem{prop}{Proposition}[section]

\newtheorem{lemma}[prop]{Lemma}
\newtheorem{cor}[thm]{Corollary}
\newtheorem{defn}[prop]{Definition}
\newtheorem*{defn*}{Definition}

\newtheorem*{claim*}{Claim}

\newcommand{\R}{{\mathbb R}}

\newcommand{\I}{{\mathcal I}}
\newcommand{\J}{{\EuScript J}}
\newcommand{\K}{{\mathcal K}}

\newcommand{\W}{{\mathcal W}}
\newcommand{\Z}{{\EuScript Z}}

\newcommand{\w}{\omega}

\newcommand{\Jbar}{\overline{J}}
\newcommand{\pibar}{\overline{\pi}}
\newcommand{\thetabar}{\overline{\theta}}

\newcommand{\Omegabar}{\overline{\Omega}}

\newcommand{\Ibar}{\overline{\I}}
\newcommand{\Mbar}{\overline{\scM}}

\newcommand{\oo}{\overline{0}}

\def\&{\wedge}

\newcommand{\di}{\partial}

\def\BB/{B\"acklund}
\def\BT/{\BB/ transformation}
\def\MA/{Monge-Amp\`ere}

\newcommand{\om}{\omega}
\newcommand{\alp}{\alpha}
\newcommand{\bet}{\beta}

\newcommand{\calC}{{\mathcal C}}
\def\Cprime#1{{\calC'_{#1}}}

\newcommand{\calK}{{\mathcal K}}

\newcommand{\calP}{{\mathcal P}}

\newcommand{\scB}{{\EuScript B}}

\newcommand{\scM}{{\EuScript M}}
\newcommand{\scN}{{\EuScript N}}

\newcommand{\scrU}{{\EuScript U}}
\newcommand{\scrV}{{\EuScript V}}

\newcommand{\B}{\scB}

\newenvironment{rem}{\begin{trivlist} \item[] {\it Remark.}}{\end{trivlist}}
\newenvironment{pf}{\begin{trivlist} \item[] {\it Proof.} }{\qed \end{trivlist} }

\begin{document}
\title{B\"acklund Transformations and Darboux Integrability for Nonlinear Wave Equations}
\author{Jeanne N. Clelland}
\address{Dept. of Mathematics, 395 UCB, University of
Colorado\\
Boulder, CO 80309-0395}
\email{Jeanne.Clelland@colorado.edu}
\author{Thomas A. Ivey}
\address{Dept. of Mathematics, College of Charleston\\
66 George St., Charleston SC 29424-0001}
\email{IveyT@cofc.edu}

\subjclass[2000]{Primary 37K35 58J72; Secondary 35L10, 37K35, 53C10, 58A15}
\keywords{B\"acklund transformations, hyperbolic Monge-Amp\`ere systems,
Weingarten surfaces, exterior differential systems, Cartan's method of equivalence}

\begin{abstract}
We prove that second-order \MA/ equations for one function of two variables
are connected to the wave equation by a \BT/ if and only if they are integrable
by the method of Darboux at second order.
\end{abstract}
\maketitle

\section{Introduction}
Roughly speaking, a \BT/ is a method for generating new solutions for a given partial
differential equation by starting with a `seed' solution to the same (or a different)
PDE and solving an auxiliary system of ODEs.  B\"acklund's original example was a
transformation that produced new pseudospherical surfaces from old, and it is equivalent
to the following system:
\begin{equation} \label{SGsingleBT}
\begin{aligned}
v_x -u_x &= \tfrac12\sin((u+v)/2), \\
v_y +u_y&= -\tfrac12\sin((u-v)/2).
\end{aligned}
\end{equation}
Given an arbitrary smooth function $u(x,y)$, this overdetermined system for $v$ is inconsistent.  However,
if $u$ satisfies the {\em sine-Gordon equation} $u_{xy}=\sin u$
then the system is consistent, and
the function $v(x,y)$, determined up to a constant of integration,
will also satisfy the sine-Gordon equation.  The transformation works in reverse, too:
given a solution $v(x,y)$ for sine-Gordon, the system determines a 1-parameter family of solutions $u(x,y)$ for the same PDE.

It is this type of B\"acklund transformation---connecting solutions of two
second-order \MA/ PDE in the plane, not necessarily the same equation---which is the general subject of
this paper.  (A second-order \MA/
equation is a PDE where the highest-order derivatives may appear nonlinearly but
only in the form of the determinant of the Hessian.)
Another important example of this type is the system
\begin{equation}\label{Liouvillewave}
\begin{aligned}
z_x -u_x&=  - 2\exp((u+z)/2), \\
z_y +u_y&= \exp((u-z)/2).
\end{aligned}
\end{equation}
In this example, if $z(x,y)$ satisfies the wave equation (in characteristic coordinates, $z_{xy}=0$),
then  the system determines a 1-parameter family of solutions of
{\em Liouville's equation} $u_{xy}=e^u$, and conversely.   \BT/s of this subtype---where one of the two
PDE involved is the standard wave equation---are the specific concern of this paper.

Liouville's equation also has the rare property that it is {\em Darboux-integrable}---in other
words, it can be solved by the method of Darboux.  (This will be defined below.)
The main point of this paper is that this is not a coincidence; more precisely, we will prove

\begin{thm}\label{ourtheorem}
Let  $(\scM^5,\I)$ be a hyperbolic \MA/ system.  If there is a {\em normal} \BT/ with 1-dimensional fibers
linking this system with the wave equation, then the first prolongation of $\I$ is
Darboux-integrable.  Conversely,  if the first prolongation of $\I$ is Darboux-integrable,
then near any point $p \in \scM$ there is an open set $U\subset\scM$ around $p$ such that
the restriction of $\I$ to $U$ is linked to the wave equation by a normal \BT/.
\end{thm}

The technical terms in this theorem must be explained.  Any single PDE or system of PDE
may be re-cast as an {\em exterior differential system (EDS)}
or differential ideal (i.e., an ideal, with respect to wedge product,
 in the ring of differential forms on a manifold, that is also closed
 under the exterior derivative),
in a way that solutions are in one-to-one correspondence
with submanifolds to which all the differential forms in the EDS pull back to be zero.  
(These submanifolds, which  must usually also satisfy a nondegeneracy condition,
are known as {\em integral manifolds} of the EDS.)
In particular, a \MA/ equation in the plane can be re-cast as the following type of EDS:

\begin{defn*}A {\em \MA/ exterior differential system} is a differential ideal $\I$
on a 5-dimensional manifold $\scM$, such that near any point of $\scM$,
$\I$ is generated algebraically by one 1-form $\theta$
 and two 2-forms $\Omega_1,\Omega_2$.
(Hence, $d\theta$ must equal a linear combination of the $\Omega$'s, plus possibly
a wedge product with $\theta$ as factor.)
The 1-form $\theta$ is required to be a {\em contact form}, i.e.,
$\theta \& d\theta \& d\theta \ne 0$.
A \MA/ system is {\em hyperbolic} if the $\Omega$'s may be chosen so that both are decomposable.
\end{defn*}

For example, for Liouville's equation we may take $\scM$ to
be $\R^5$ with coordinates $x,y,u,p,q$, and let
\begin{equation}\label{Liouvilleforms}
\theta=du-p\,dx-q\,dy,\qquad
\Omega_1=(dp - e^u dy) \& dx,\qquad
\Omega_2=(dq - e^u dx) \& dy.
\end{equation}
(Note that $d\theta = -\Omega_1 - \Omega_2$.)
Given a solution $u=f(x,y)$ of the PDE, we can construct
a surface $\Sigma \subset \R^5$ such that $i^*\theta=0$, $i^*\Omega_1 = i^*\Omega_2 = 0$
(where $i:\Sigma \hookrightarrow \R^5$ is the inclusion map)
by setting $u=f(x,y)$, $p=f_x(x,y)$ and $q=f_y(x,y)$.  Conversely,
any surface $\Sigma$ satisfying $i^*\theta=0$, $i^*\Omega_1 =i^*\Omega_2=0$ and
the nondegeneracy condition $i^* dx \& dy \ne 0$
is the graph of a solution constructed in this way.

In the body of the paper, we will also use another type of EDS:
\begin{defn*} A {\em Pfaffian exterior differential system} is a differential
ideal $\I$ on an arbitrary manifold $\scM$, defined by a vector bundle $I \subset T^*\scM$,
such that a differential form belongs to $\I$ if and only if it is a linear combination
of wedge products involving sections of $I$ or their exterior derivatives.
(In practice, our Pfaffian systems will be specified by giving a list of
1-forms that span the fiber of $I$ at each point.)  The {\em rank} of a Pfaffian
system is the rank of the vector bundle.
\end{defn*}

A Pfaffian system satisfies the {\em Frobenius condition} or is said to
be {\em integrable} if the exterior derivative
of any section of $I$ is in the {\em algebraic} ideal generated
by $I$.
Any Frobenius system is locally equivalent to a (possibly underdetermined)
system of ordinary differential equations; see Chapter 1 in \cite{CFB}.

\bigskip
Theorem \ref{ourtheorem} is about relations between exterior differential systems; in particular,
we have the following definition from \cite{CFB}:
\begin{defn}\label{general-Backlund-defn}
A {\em \BT/ between two exterior differential systems} $(\scM,\I)$ and
$(\Mbar,\Ibar)$ is a manifold $\B$ equipped with submersions $\pi:\B\to \scM$
and $\pibar:\B \to \Mbar$ (see diagram below) and vector bundles $J,\Jbar \subset T^*\B$
such that
\newline
(i) the fibers of $\pi$ and $\pibar$ are transverse in $\B$;
\newline
(ii) the rank of $J$ equals the dimension of the fibers of $\pi$, and
sections of $J$ pull back to the fibers of $\pi$ to span the
cotangent space of each fiber;
\newline
(iii) $\Jbar$ is similarly related to the fibers of $\pibar$;
\newline
(iv) the algebraic ideal $\J$ generated by $\pi^* \I$ and sections of $J$ is
the same as the algebraic ideal generated by $\pibar^* \Ibar$ and sections of $\Jbar$,
and $\J$ is a differential ideal.
\end{defn}

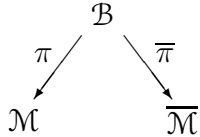
\setlength{\unitlength}{2pt}
\begin{center}
\begin{picture}(40,30)(0,0)
\put(5,5){\makebox(0,0){$\scM$}}
\put(35,5){\makebox(0,0){$\Mbar$}}
\put(20,25){\makebox(0,0){$\B$}}
\put(16,21){\vector(-3,-4){9}}
\put(24,21){\vector(3,-4){9}}
\put(7,16){$\pi$}
\put(30,16){$\pibar$}
\end{picture}
\end{center}

The impact of the last condition is that if $\scN \subset \scM$
is an integral manifold of $\I$, then sections of $J$ pull back to
$\pi^{-1}(\scN)$ to satisfy the Frobenius condition, so that
integral manifolds of $\J$ inside $\pi^{-1}(\scN)$ may be constructed by
solving ODE; moreover, the image under $\pibar$ of each of 
these integral manifolds is an integral manifold
of $\Ibar$.  Because the definition is symmetric, this also works in the other direction: 
given an integral manifold of $\Ibar$, 
we can solve a Frobenius system on the inverse image in $\B$
to obtain a family of integral manifolds of $\I$.
For example, given a solution $z(x,y)$ of the wave equation,
substitution in \eqref{Liouvillewave} gives an overdetermined system of ODE for
a solution $u(x,y)$ of Liouville's equation, and in this context the
Frobenius condition is exactly the compatibility condition for the ODE system.

The condition of {\em normality} for \BT/s, assumed in Theorem \ref{ourtheorem},
will be explained in \S2.

\bigskip

A hyperbolic \MA/ system is a special case of {\em hyperbolic EDS}:
\begin{defn*}A {\em hyperbolic EDS of class $k$} is a differential ideal defined on
a manifold of dimension $k+4$ that,\ near any point of manifold,
is generated algebraically by $k$ 1-forms and two decomposable 2-forms.
\end{defn*}

Associated to a given hyperbolic EDS $\I$ of class $k$ are two characteristic distributions,
one corresponding to each decomposable 2-form generator.
At each point, the distribution is given by the 2-dimensional subspace of the
tangent space annihilated by the $k$ 1-forms of the system
and the factors of the chosen decomposable 2-form.  (These annihilators
form a rank $k+2$ Pfaffian system, known as a {\em characteristic system} of $\I$.)
A hyperbolic EDS $\I$ is {\em integrable by the method of Darboux},
or {\em Darboux-integrable} for short, if both characteristic distributions have two
independent first integrals, i.e., functions which are constant
along all curves tangent to the distribution, and whose differentials are
pointwise linearly independent from the 1-forms of $\I$.
Such functions are also known as {\em characteristic invariants}, since they are constant
along integral curves of the distribution.

The Darboux-integrability condition has the virtue that it is easy to check, using
only differentiation and linear algebra, by calculating the successive derived
systems of each characteristic system.  An extensive discussion
of hyperbolic EDS and Darboux-integrability, with worked-out examples, is
available in Chapter 6 of \cite{CFB}.
For the purposes of this paper, we will need a few more facts about the method of Darboux:
\begin{itemize}
\item It is known that any \MA/ system which is Darboux-integrable (i.e.,
has two characteristic invariants for each distribution) is equivalent
to the standard wave equation under a contact transformation (see, e.g., Thm. 2.1
in \cite{BGG}).
\item If a hyperbolic \MA/ system $\I$ has a pair of independent first integrals for exactly
one of its characteristic distributions, then $\I$ is said to be {\em integrable
by the method of Monge} or {\em Monge-integrable} for short.
(The analogous term for hyperbolic EDS of arbitrary class $k$ is {\em Darboux semi-integrable}.)
\item If a hyperbolic EDS of class $k$ fails to be Darboux-integrable, it is possible that
its {\em prolongation}, which is a hyperbolic EDS of class $k+2$, is Darboux-integrable.
Thus, a given \MA/ system may lead to a hyperbolic system that is Darboux-integrable only after
sufficiently many prolongations.
\end{itemize}
Prolongation of an EDS is essentially the process
of adding higher derivatives as new variables and adjoining to the ideal
the differential equations satisfied by the higher derivatives.
For example,
for Liouville's equation we add variables $r$ and $t$ to stand for $u_{xx}$ and $u_{yy}$
respectively, and adjoin the 1-forms $\theta_1 = dp - r\,dx -e^u dy$,
$\theta_2 = dq - e^u dx - t\,dy$.  The new ideal is a Pfaffian system on $\R^7$
(with coordinates $x,y,u,p,q,r,t$) generated by 1-forms $\theta_0,\theta_1,\theta_2$.
(Note that
the 2-forms $\Omega_1$ and $\Omega_2$ given in \eqref{Liouvilleforms} are now
in the ideal generated algebraically by these $\theta_0,\theta_1,\theta_2$.)
The set of algebraic generators of the new ideal is
completed by computing the exterior derivatives of $\theta_1,\theta_2$,
and these are expressible as linear combinations
of the decomposable forms $\Omega'_1 = (dr - p e^u dy) \& dx$, $\Omega'_2 =
(dt - q e^u dx) \& dy$, modulo multiples of $\theta_0,\theta_1,\theta_2$.
(Thus, the new ideal is a hyperbolic EDS of class 3.)
Let $\Delta_1, \Delta_2$ be the corresponding characteristic distributions for
the prolongation (i.e., $\Delta_i$ is annihilated by $\theta_0,\theta_1,\theta_2$
and the factors of $\Omega'_i$).
To see that the system is Darboux-integrable,
note that $x, r - \tfrac12 p^2$ are invariants for $\Delta_1$ and
$y, t -\tfrac12 q^2$ are invariants for $\Delta_2$.

\begin{rem}  Both Darboux-integrability and the transformation \eqref{Liouvillewave} enable
one to express all solutions of Liouville's equation via specifying two arbitrary functions
and integrating systems of ODE, and these two solution methods are
equivalent, although Darboux's method requires one to solve more ODEs.
For, as mentioned above, substituting the wave equation solution $z=f(x)+g(y)$
in \eqref{Liouvillewave} produces two
compatible ODEs for $u(x,y)$.
Given an initial value for $u$, these
may be integrated simultaneously in the $x$- and $y$-directions
to propagate a solution over an open set in the $xy$-plane.
On the other hand, under Darboux's method we obtain ODEs by setting one invariant
in each characteristic system to be an arbitrary function of the other, yielding in this case
the equations
$$p_x-\tfrac12 p^2 = \phi(x), \qquad q_y-\tfrac12 q^2 = \psi(y),$$
which, together with $u_x=p$ and $u_y=q$, may be integrated to obtain the solution.
(In fact, the data for these two methods are related by $\phi=f''$ and $\psi=g''$, but
in other cases we cannot expect the relationship to be this simple.)
\end{rem}

Next, we will put Theorem \ref{ourtheorem} in context with other results both classical and modern.
Much of what was known in the 19th century about solving second-order PDE for
one function of two variables was summarized in Goursat's treatise \cite{GLecons}.
In Volume 2, \S181 of that work, we find the following result:
\begin{thm}[Darboux-Goursat]
Suppose that a second-order PDE for $z$ as a function of $x,y$ has the property
that there exists a Pfaffian system
$$dF_i = \Phi_i\,d\alpha + \Psi_i d\beta, \qquad 1\le i \le \ell,$$
and formulas
\begin{equation}\label{dgxv}
x = V_1, \quad y = V_2, \quad z=V_3,
\end{equation}
where $\Phi_i$, $\Psi_i$, $V_1, V_2, V_3$ are functions of $F_1,\ldots, F_\ell,$ $\alpha,\beta$,
$f(\alpha)$, $g(\beta)$ and finitely many of the derivatives of $f$ and $g$,
such that \eqref{dgxv} satisfies the Frobenius condition
for arbitrary choices of functions $f$ and $g$, and gives an implicit solution of the PDE for
arbitrary choices of initial data for the Frobenius system.
Then the PDE is Darboux-integrable after finitely many prolongations.
\end{thm}
The hypotheses of the Darboux-Goursat theorem are fulfilled if the given PDE is
linked to the standard wave equation by a \BT/.  (For, the d'Alembert formula
gives solutions of the wave equation $u_{\alpha\beta}=0$ in the form
$u=f(\alpha)+g(\beta)$ for arbitrary $f$ and $g$, and the Pfaffian
system in the theorem is given by the equations of the \BT/.)
Compared with the Darboux-Goursat theorem, one direction of our theorem has a stronger hypothesis
(essentially, that $\ell=2$ and only first derivatives of $f$ and $g$ are involved)
and a stronger conclusion (that at most one prolongation is required to get Darboux-integrability).

In Theorem 6.5.14 in \cite{CFB} it is shown, by an elementary argument, that
Darboux-integrability of the prolongation implies that there is a \BT/ between the
prolongation (not the original \MA/ system, but one defined by an EDS on a 7-dimensional
manifold) and the wave equation (defined by an EDS on a 5-dimensional manifold).
However, this asymmetric transformation--relating the 2-jets of solutions of
one PDE to the 1-jets of another---is less than satisfying, compared to
more symmetrical transformations like \eqref{Liouvillewave}.  Our analysis in \S4 shows
that it is a much more delicate matter to show that there exists a \BT/ between
two \MA/ systems.  We should also note that
the argument given in \cite{CFB} for the other direction (\BB/-equivalence
to the wave equation implies Darboux-integrability) is unfortunately incorrect,
and the proof we give in \S3 of this paper is along completely different lines.

\bigskip

We now briefly outline the rest of the paper.  In \S2 we set up the basic machinery required for
the first half of the proof, namely, the $G$-structure for \BT/s originally introduced by the first author
in \cite{C1}.  In \S3 we prove the forward direction in our theorem by
following the implications (for the invariants of the $G$-structure) of the existence of a \BT/ to the wave equation.
In \S4 we prove the converse direction by constructing, for any given Darboux-integrable Monge-Amp\`ere equation,
an involutive exterior differential system whose solutions are such transformations.
In \S5, we discuss our results in the context of earlier classifications
of Darboux-integrable equations and of \BT/s to the wave equation; 
we also outline an alternate
proof technique for the converse direction, which can in some examples be used
to establish global existence of the transformation.
In \S6 we discuss further steps in our research program.

We are grateful to the referee who read the first version of this paper, and gave us many 
useful comments and suggestions.

\section{$G$-structure for \BT/s}\label{wsetup}

The material in this section is taken from the first author's paper \cite{C1}; additional details may be found there.

Suppose that $(\scM, \I)$ and $(\Mbar, \Ibar)$ are hyperbolic
Monge-Amp\`ere systems, with
\[ \I = \{ \theta, \Omega_1, \Omega_2 \}, \qquad \Ibar = \{ \thetabar, \Omegabar_1, \Omegabar_2 \}.\]
As a special case of Definition \ref{general-Backlund-defn}, we define a B\"acklund transformation between $(\scM, \I)$ and $(\Mbar,
\Ibar)$ to be a 6-dimensional submanifold $\B \subset \scM \times \Mbar$
for which the pullbacks to $\B$ of the forms
$\Omega_1, \Omega_2, \Omegabar_1, \Omegabar_2$ have the property that
\begin{equation*}
 \Omega_i \equiv \Omegabar_i \mod{\{\theta, \thetabar\}}, \ \ i=1,2.
\end{equation*}
(The vector bundles $J,\Jbar \subset T^*\scB$ mentioned in Definition \ref{general-Backlund-defn}
are in this case spanned by the pullbacks of $\thetabar$ and $\theta$, respectively.)
A \BT/ is {\em normal} in the sense of Theorem \ref{ourtheorem} if
the pullbacks to $\B$ of the 2-forms $d\theta$, $d\thetabar$ are linearly independent modulo $\{\theta, \thetabar\}$.

Now let $\J$ be the ideal on $\B$ generated by the pullbacks of $\I$ and $\Ibar$; according to the conditions above, $\J$ is generated algebraically by the forms $\{\theta, \thetabar, d\theta, d\thetabar \}. $

Since $\I$ and $\Ibar$ are hyperbolic, locally there exist 1-forms
$\om^1, \om^2,
\om^3, \om^4$ on $\B$ such that $\{\theta,\, \thetabar,\, \om^1,\,
\om^2,\, \om^3,\, \om^4\}$ is a coframing of $\B$ (i.e., a set of
1-forms that restricts, at each point, to be a basis for the cotangent
space of $\B$)  and
\[ \J = \{\theta,\, \thetabar,\, \om^1 \& \om^2,\, \om^3 \& \om^4\}.
\]
(It is important to note that $\theta$ and $\thetabar$ are each
{\em separately} determined up to a scalar multiple, since each
determines the contact structure on a 5-manifold.)
Any such coframing has the property that
\begin{align*}
d\theta & \equiv A_1\, \om^1 \& \om^2 + A_2\,\om^3 \& \om^4
\mod{\{\theta,\, \thetabar\}},\\
d\thetabar & \equiv A_3\, \om^1 \& \om^2 + A_4\, \om^3 \& \om^4
\mod{\{\theta,\, \thetabar\}}
\end{align*}
for some nonvanishing functions $A_1, A_2, A_3, A_4$.  Since $d\theta,\,
d\thetabar$ are required to be linearly independent 2-forms at each
point of $\B$, we must have $A_1 A_4 - A_2 A_3 \neq 0$.

By rescaling the $\om^i$ and adding multiples of $\theta$ and
$\thetabar$ to the $\om^i$ if necessary, we can arrange that
\begin{align}
d\theta & \equiv A_1\, \om^1 \& \om^2 + \om^3 \& \om^4
\mod{\theta}, \label{dtheta}\\
d\thetabar & \equiv \om^1 \& \om^2 + A_2\, \om^3 \& \om^4
\mod{\thetabar} \notag
\end{align}
for some nonvanishing functions $A_1,\, A_2$ on $\B$ with $A_1
A_2 \neq 1$.
This coframing is not unique; any other such coframing
$\{\tilde{\theta},\, \tilde{\thetabar},\, \tilde{\om}^1,\,
\tilde{\om}^2,\,
\tilde{\om}^3,\, \tilde{\om}^4\}$ has the form
\begin{equation}
 \begin{bmatrix} \tilde{\theta}\\[0.1in] \tilde{\thetabar}\\[0.1in]
\tilde{\om}^1\\[0.1in] \tilde{\om}^2\\[0.1in] \tilde{\om}^3\\[0.1in]
\tilde{\om}^4
\end{bmatrix} =   \begin{bmatrix} b_{11}b_{22} - b_{12}b_{21} & 0 &
0 & 0 & 0 & 0 \\[0.1in] 0 &  a_{11}a_{22} - a_{12}a_{21} & 0 & 0 & 0 & 0
\\[0.1in] 0 & 0 & a_{11} & a_{12} & 0 & 0 \\[0.1in] 0 & 0 & a_{21} &
a_{22} & 0 & 0
\\[0.1in] 0 & 0 & 0 & 0 & b_{11}&  b_{12} \\[0.1in] 0 & 0 & 0 &
0 & b_{21} & b_{22}
\end{bmatrix}^{-1} \begin{bmatrix} \theta\\[0.1in] \thetabar\\[0.1in]
\om^1\\[0.1in]
\om^2\\[0.1in] \om^3\\[0.1in] \om^4 \end{bmatrix},  \label{G0freedom}
\end{equation}
where
$ b_{11}b_{22} - b_{12}b_{21} \neq 0, \
a_{11}a_{22} - a_{12}a_{21} \neq 0.$  (The inverse is included for
greater ease of computation in carrying out the method of equivalence.)  A coframing satisfying
\eqref{dtheta} is called {\em adapted}, and the group $G$ of matrices
of the
above form is called the {\em structure group} of the equivalence
problem.  (In fact, the most general choice of structure group would
include a discrete component interchanging the
distributions $\{\om^1, \om^2\}$ and $\{\om^3, \om^4\}$.  However,
this freedom does not contribute anything crucial to the structure
group, and it is easier to work with a connected group.)  The associated {\em $G$-structure} is the principal $G$-bundle $\calP \to \B$  whose local sections are precisely the adapted coframings over a neighborhood of $\B$.

In \cite{C1}, it is shown that $\calP$ has {\em structure equations}
\begin{multline}
\begin{bmatrix} d\theta\\[0.1in] d\thetabar\\[0.1in] d\om^1\\[0.1in]
d\om^2\\[0.1in] d\om^3\\[0.1in] d\om^4
\end{bmatrix} =
-\begin{bmatrix}\bet_1 + \bet_4 & 0 & 0 & 0 & 0 & 0 \\[0.1in] 0 & \alp_1 +
\alp_4 & 0 & 0 & 0 & 0 \\[0.1in] 0 & 0 & \alp_1 & \alp_2 & 0 & 0
\\[0.1in] 0 & 0 &
\alp_3 & \alp_4 & 0 & 0 \\[0.1in] 0 & 0 & 0 & 0 & \bet_1 & \bet_2
\\[0.1in] 0 & 0 & 0 &
0 & \bet_3 & \bet_4 \end{bmatrix} \&
\begin{bmatrix} \theta\\[0.1in] \thetabar\\[0.1in] \om^1\\[0.1in]
\om^2\\[0.1in] \om^3\\[0.1in] \om^4
\end{bmatrix}  \\[0.2in] \label{G0struct}
+ \begin{bmatrix} \theta \& (A_1 C_2\, \om^1 - A_1 C_1\, \om^2) + A_1\, \om^1 \& \om^2 + \om^3 \& \om^4
\\[0.1in]
\thetabar \& (A_2 C_4\, \om^3 - A_2 C_3\, \om^4) + \om^1 \& \om^2 + A_2\, \om^3 \&
\om^4 \\[0.1in]
B_1\, \theta \& \thetabar + C_1\, \om^3 \& \om^4  \\[0.1in]
B_2\, \theta \& \thetabar + C_2\, \om^3 \& \om^4  \\[0.1in]
B_3\, \theta \& \thetabar + C_3\, \om^1 \& \om^2  \\[0.1in]
B_4\, \theta \& \thetabar + C_4\, \om^1 \& \om^2 \end{bmatrix}
\end{multline}
for some functions $A_i, B_i, C_i$ and 1-forms  $\alp_i, \bet_i$ on $\calP$.
These equations are chosen so
that the matrix in \eqref{G0struct} takes values in the Lie algebra
$\mathfrak{g}$
of $G$; this is a standard step in the method of equivalence.  (See
\cite{G89} for details.)

The 1-forms $\alp_i, \bet_i$ are linearly independent from each other and from $\theta, \thetabar, \om^i$; they are called {\em pseudoconnection forms,} or more concisely (but imprecisely) {\em connection forms} on $\calP$.  They are well-defined only up to transformations of the form
\begin{alignat}{2}
\alp_1 & \mapsto \alp_1 + r_1\, \om^1 + r_2\, \om^2, \qquad \qquad & \bet_1 &
\mapsto \bet_1 + s_1\, \om^3 + s_2\, \om^4, \notag \\
\alp_2 & \mapsto \alp_2 + r_2\, \om^1 + r_3\, \om^2, \qquad \qquad & \bet_2 &
\mapsto \bet_2 + s_2\, \om^3 + s_3\, \om^4, \label{G0connfreedom} \\
\alp_3 & \mapsto \alp_3 + r_4\, \om^1 - r_1\, \om^2, \qquad \qquad & \bet_3 &
\mapsto \bet_3 + s_4\, \om^3 - s_1\, \om^4, \notag \\
\alp_4 & \mapsto \alp_4 - r_1\, \om^1 - r_2\, \om^2, \qquad \qquad & \bet_4 &
\mapsto \bet_4 - s_1\, \om^3 - s_2\, \om^4.\notag
\end{alignat}

\begin{rem}
The coefficients $A_i, B_i, C_i$ are called {\em torsion functions}.  They may be
interpreted as the components of well-defined tensors associated to the \BT/, as follows.

A hyperbolic \MA/ system naturally equips the underlying manifold
$\scM^5$ with a line bundle $L \subset T^*\scM$ and
two rank 3 characteristic bundles $K_1, K_2 \subset T^*\scM$ whose intersection is $L$.
(The generator 1-form $\theta$ is a section of $L$, and the factors of the
decomposable generator 2-forms $\Omega_1$ and $\Omega_2$ span a complement of $L$ within $K_1$ and $K_2$,
respectively.)  The $G$-structure for the \BT/ shows that $\scB^6$ is equipped with
a well-defined splitting of its cotangent bundle:
\def\Lbar{\overline{L}}
\begin{equation}\label{TBsplit}
T^*\scB = L \oplus \Lbar \oplus W_1 \oplus W_2,
\end{equation}
where $L$ and $\Lbar$ are the pullbacks of the \MA/ line bundles from $\scM$ and $\Mbar$
respectively, and $W_1, W_2$ are spanned by $\w^1, \w^2$ and $\w^3, \w^4$ respectively.
The normal \BB/ condition implies that
$$\pi^* K_i = L \oplus W_i, \qquad \pibar^* \overline{K}_i = \Lbar \oplus W_i, \qquad i=1,2.$$

The splitting \eqref{TBsplit} induces a corresponding splitting of $\Lambda^2 T^* \scB$,
whose summands include the line bundles $L \otimes \Lbar$, $\Lambda^2 W_1$ and $\Lambda^2 W_2$.  We may then define
a natural map from sections of $L$ to sections of $\Lambda^2 W_1$, given by
$$\delta_0: \theta \mapsto \text{projection of } d\theta \text{ into } \Lambda^2 W_1.$$
But this map is linear under multiplication by functions, and so gives a well-defined
map between the corresponding vector bundles.
The structure equations \eqref{G0struct}
show that $A_1 \w^1 \& \w^2$ is the value of $\delta_0$ applied to the first
member $\theta$ of the coframe.  Hence, $A_1$ is the component, with respect to
the give coframe, of a well-defined tensor in $L^* \otimes \Lambda^2 W_1$.
Similarly, $A_2$ is the component of a well-defined tensor
in $\Lbar^* \otimes \Lambda^2 W_2$

We may similarly define a map on sections of $W_1$ by
$$\delta_1: \w \mapsto \text{projection of } d\w \text{ into } L_1 \otimes L_2,$$
which again is linear under multiplication by functions.  The structure equations
show that $B_1 \theta \wedge \thetabar$ and $B_2 \theta \wedge \thetabar$ give
the value of $\delta_1$ on the basis sections $\w^1, \w^2$ of $W_1$ respectively.
Thus, the vector $[B_1\ \ B_2]$ gives the components, with respect to the coframe, of a well-defined
tensor in $W_1^* \otimes L_1 \otimes L_2$.  In a similar way, we see that
$[B_3\ \ B_4]$ are the components of a tensor in $W_2^* \otimes L_1 \otimes L_2$,
$[C_1\ \ C_2]$ are the components of a tensor in $W_1^* \otimes \Lambda^2 W_2$,
and $[C_3\ \ C_4]$ are the components of a tensor in $W_2^* \otimes \Lambda^2 W_1$,
all defined by taking the exterior derivative of a section and projecting
into the appropriate summand of $\Lambda^2 T^* \scB$.
\end{rem}
\bigskip

The following results are proved in \cite{C1}:

\begin{prop}
The vectors $[B_1 \ \ B_2], \ [B_3 \ \ B_4], \ [C_1 \ \ C_2], \ [C_3 \ \ C_4]$ are {\em relative invariants:} given any point $m \in \B$, they are each either zero for every adapted coframing at $m$, or nonzero for every adapted coframing at $m$.
\end{prop}

\begin{prop}
If $[C_1 \ \ C_2] = [C_3 \ \ C_4] = [0\ \ 0]$, then $[B_1 \ \ B_2] = [B_3 \ \ B_4] = [0\ \ 0]$ as well, and $(\scM, \I)$ and $(\Mbar, \Ibar)$ are each contact equivalent to the Monge-Amp\`ere system representing the standard wave equation.
\end{prop}

\begin{prop}\label{one-C-vector-zero}
If $[C_1 \ \ C_2] = [0\ \ 0]$ (resp., $[C_3 \ \ C_4] = [0\ \ 0]$), then $[B_1 \ \ B_2] = [0\ \ 0]$ (resp., $[B_3 \ \ B_4] = [0\ \ 0]$) as well, and $(\scM, \I)$ and $(\Mbar, \Ibar)$ are each Monge-integrable.
\end{prop}

\begin{prop}
Suppose that the vectors $[C_1 \ \ C_2]$ and $[C_3 \ \ C_4]$ are both nonzero.  Then the vectors $[B_1 \ \ B_2]$ and $[B_3 \ \ B_4]$ are either both zero or both nonzero.
\end{prop}

If $[B_1 \ \ B_2] = [B_3 \ \ B_4] = [0\ \ 0]$, then the differential ideal generated by $\{\om^1, \om^2, \om^3, \om^4\}$ is a Frobenius system.  (The converse is true as well.) It follows that locally, there exists a 4-manifold $\scrV$ which is a quotient of
$\B$ and for which the 1-forms $\om^1, \om^2, \om^3, \om^4$ are
semi-basic for the projection $\rho:\B \to \scrV$.  (Here ``locally'' refers
to the fact that any point in $\B$ has a neighborhood which
possesses such
a quotient, and ``semi-basic'' means that the
restrictions of the $\om^i$ to the fibers of the projection vanish
identically.  See \cite{CFB} for details.)  In fact, this quotient
factors through each of the quotients $\pi: \B \to \scM$ and $\pibar: \B \to \Mbar$, as shown
by the following commutative diagram.
\setlength{\unitlength}{2pt}
\begin{center}
\begin{picture}(40,40)(0,10)
\put(5,28){\makebox(0,0){$\scM$}}
\put(35,28){\makebox(0,0){$\Mbar$}}
\put(20,45){\makebox(0,0){$\B$}}
\put(16,41){\vector(-1,-1){9}}
\put(24,41){\vector(1,-1){9}}
\put(7,37){$\pi$}
\put(31,37){$\pibar$}
\put(20,40){\vector(0,-1){23}}
\put(23,27){\makebox(0,0){$\rho$}}
\put(20,13){\makebox(0,0){$\scrV$}}
\put(7,25){\vector(1,-1){9}}
\put(33,25){\vector(-1,-1){9}}
\end{picture}
\end{center}
The vanishing of the vector $[B_1\ \ B_2]$ implies that the
span of $\{\w^1, \w^2\}$ is unchanged along the fibers of $\rho$,
and is thus the pullback via $\rho$ of a well-defined rank 2 sub-bundle of $T^*\scrV$.
When the vector $[B_3 \ \ B_4]$ also vanishes,
the ideal $\{\om^1 \& \om^2,\  \om^3 \& \om^4\}$
is the pullback via $\rho$ of a well-defined hyperbolic system
$\mathcal{H}$ of class 0 on $\scrV$,
and $\I, \Ibar$ are both integrable extensions of $\mathcal{H}$.

\BT/s of this type are called {\em holonomic}.  One can test whether a \BT/ is holonomic by checking whether the Pfaffian system on $\B$ spanned by the intersection of the basic forms for $\pi$ with the basic forms for $\pibar$ is Frobenius.  Note that
the basic  forms for $\pi$ are spanned by the {\em Cartan system}\footnote{The {\em Cartan system} for a given EDS $\I$ is the smallest Frobenius system that contains $\I$.} of $\I$, and the basic forms for $\pibar$ are spanned by the Cartan system of $\Ibar$.

Holonomic \BT/s are generally considered less interesting than non-holonomic \BT/s because of their limited capacity to generate new solutions,
which we now explain. Given an integral surface $\scN$ of $(\scM, \I)$, solving the Frobenius system $\J$ on $\pi^{-1}(\scN)$ produces a 1-parameter family of integral surfaces $\overline{\scN}_{\lambda}$ of $(\Mbar, \Ibar)$.  Reversing the process, beginning with any one of the integral manifolds $\overline{\scN}_{\lambda}$, in turn produces a 1-parameter family of integral surfaces $\scN_{\lambda, \mu}$ of $(\scM, \I)$.  In general, this results in a 2-parameter family of integral surfaces of $(\scM, \I)$, and iterating the process produces an ever-increasing family of new integral surfaces for both Monge-Amp\`ere systems.  

For example, consider the system \eqref{Liouvillewave}.  If we substitute the trivial solution $z(x,y)=0$ of the wave equation into \eqref{Liouvillewave}, the resulting overdetermined PDE system for $u$ yields the 1-parameter family of solutions
\begin{equation}
u(x,y) = -2\ln (-x -\tfrac{1}{2} y - c_1) \label{Liouville-u1}
\end{equation}
to Liouville's equation.  Reversing the process, substituting \eqref{Liouville-u1} into \eqref{Liouvillewave} produces a PDE system for $z$ which has a 2-parameter family of solutions:
\begin{equation}
z(x,y) = 2 \ln (-y - c_2) - 2 \ln (2x + 2c_1 - c_2) . \label{Liouville-z2}
\end{equation}
Finally, substituting \eqref{Liouville-z2} into \eqref{Liouvillewave} produces a PDE system for $u$ which has a 3-parameter family of solutions:
\begin{equation}
u(x,y) = -2\ln \big{(}c_3 xy  + (c_2 c_3 - 1) x +
(c_1 c_1 - \tfrac{1}{2} c_2 c_3 -\tfrac{1}{2}) y +
(c_1 c_2 c_3 - \tfrac{1}{2} c_2^2 c_3 -  c_1) \big{)}. \label{Liouville-u3}
\end{equation}
It is clear that the solutions \eqref{Liouville-u1} form a proper subset of the solutions \eqref{Liouville-u3}, since the argument of the latter contains an $xy$ term.

For the system \eqref{Liouvillewave}, and for non-holonomic \BT/s in general, successive iterations of this process continue to produce new solutions. However, if the \BT/ is holonomic, then all integral surfaces of $(\scM, \I)$ and $(\Mbar, \Ibar)$ produced by this process must lie in the inverse image of a single integral surface of $(\scrV, \mathcal{H})$.  It follows that successive iterations can produce no more than a 1-parameter family of integral surfaces for each Monge-Amp\`ere system.

\section{Proof that B\"acklund implies Darboux}\label{wprovedarboux}

Now suppose that we have a \BT/ as in \S \ref{wsetup}, and that the Monge-Amp\`ere system $(\scM^5, \I)$ is contact equivalent to the standard wave equation $Z_{XY}=0$.  We can choose local coordinates $(X,Y,Z,P,Q)$ on $\scM$ such that $\I$ is generated by the forms
\[ \theta = dZ - P\, dX - Q\, dY, \ \ \ \ \ 
\Omega_1 = dX \& dP, \ \ \ \ \ 
\Omega_2 = dY \& dQ. \]
There is a unique local section $\sigma = (\theta, \thetabar, \om^1, \om^2, \om^3, \om^4):\B \to \calP$ satisfying
\begin{align}
\theta & = dZ - P\, dX - Q\, dY, \notag \\
\om^1 & = dX + C_1 \theta, \notag \\
\om^2 & = dP + C_2 \theta, \label{wavecoframe} \\
\om^3 & = dY + C_3 \theta, \notag \\
\om^4 & = dQ + C_4 \theta \notag
\end{align}
for some functions $C_i$ on $\B$.  (Note that, because specifying this portion of the coframing determines a unique local section of $\calP$, the 1-form $\thetabar$  is also uniquely determined.)  The functions $C_i$ are the pullbacks under $\sigma$ of the corresponding torsion functions, and this coframing has $A_1 \equiv 1$.  When the structure equations \eqref{G0struct} are pulled back to $\B$ via $\sigma$, the 1-forms $\alp_i, \bet_i$---which were linearly independent from the 1-forms $\theta, \thetabar, \om^i$ on $\calP$---must pull back to some linear combinations of these 1-forms.

Now we embark on the process of comparing the structure equations \eqref{G0struct} to those for the explicit coframing above.  First, note that
\[  0 = d(dX) = d(\om^1 - C_1 \theta) \equiv -(\alp_1 \& \om^1 + \alp_2 \& \om^2 + C_1 \om^1 \& \om^2)  \mod{\theta}. \]
Therefore, by choosing $r_2, r_3$ appropriately in \eqref{G0connfreedom}, we may assume that
\[ \alp_1 = a_{10}\theta + a_{11}\om^1 + \tfrac{1}{2} C_1 \om^2, \qquad
\alp_2 = a_{20}\theta - \tfrac{1}{2} C_1 \om^1 \]
for some functions $a_{10}, a_{11}, a_{20}$ on $\B$.
Similar considerations of $d(dP), d(dY), d(dQ)$ modulo $\theta$ yield similar expressions for the remaining connection forms:
\begin{alignat*}{2}
\alp_1 &= a_{10}\theta + a_{11}\om^1 + \tfrac{1}{2} C_1 \om^2,  \qquad  \qquad & \bet_1 &= b_{10}\theta + b_{13}\om^3 + \tfrac{1}{2} C_3 \om^4, \\
\alp_2 &= a_{20}\theta - \tfrac{1}{2} C_1 \om^1, \qquad \qquad & \bet_2 &= b_{20}\theta - \tfrac{1}{2} C_3 \om^3,\\
\alp_3 &= a_{30}\theta + \tfrac{1}{2} C_2 \om^2, \qquad \qquad & \bet_3 &= b_{30}\theta + \tfrac{1}{2} C_4 \om^4,\\
\alp_4 &= a_{40}\theta - \tfrac{1}{2} C_2 \om^1 + a_{42}\om^2,  \qquad \qquad & \bet_4 &= b_{40}\theta - \tfrac{1}{2} C_4 \om^3 + b_{44}\om^4 .
\end{alignat*}

Next, a straightforward computation shows that for the coframing \eqref{wavecoframe},
\[ d\theta = \theta \& (C_2 \om^1 - C_1 \om^2 + C_4 \om^3 - C_3 \om^4)
+ \om^1 \& \om^2 + \om^3 \& \om^4. \]
Comparing this with the first structure equation in \eqref{G0struct} yields
\[ b_{13} = \tfrac{3}{2} C_4 , \qquad
b_{44} = -\tfrac{3}{2} C_3.  \]
In order to continue this comparison, we need to introduce derivatives of the functions $A_2, C_i$.  So, set
\begin{align}
dA_2 & = A_{2,0}\theta + A_{2,\overline{0}}\thetabar + A_{2,1}\om^1 + A_{2,2}\om^2 + A_{2,3}\om^3 + A_{2,4}\om^4, \notag \\
dC_1 & = C_{1,0}\theta + C_{1,\overline{0}}\thetabar + C_{1,1}\om^1 + C_{1,2}\om^2 + C_{1,3}\om^3 + C_{1,4}\om^4, \notag \\
dC_2 & = C_{2,0}\theta + C_{2,\overline{0}}\thetabar + C_{2,1}\om^1 + C_{2,2}\om^2 + C_{2,3}\om^3 + C_{2,4}\om^4, \label{first-derivs} \\
dC_3 & = C_{3,0}\theta + C_{3,\overline{0}}\thetabar + C_{3,1}\om^1 + C_{3,2}\om^2 + C_{3,3}\om^3 + C_{3,4}\om^4, \notag \\
dC_4 & = C_{4,0}\theta + C_{4,\overline{0}}\thetabar + C_{4,1}\om^1 + C_{4,2}\om^2 + C_{4,3}\om^3 + C_{4,4}\om^4 . \notag
\end{align}
Comparing the structure equations \eqref{G0struct} for $d\om^i$ with the derivatives of the explicit forms $\om^i$ in \eqref{wavecoframe} yields
\begin{alignat}{2}
a_{10} & = C_{1,1} - C_1 C_2, \qquad \qquad & b_{10} & = C_{3,3} - C_3 C_4 \notag \\
a_{20} & = C_{1,2} + C_1^2, \qquad \qquad & b_{20} & = C_{3,4} + C_3^2, \notag \\
a_{30} & = C_{2,1} - C_2^2, \qquad \qquad & b_{30} & = C_{4,3} - C_4^2, \notag \\
a_{40} & = C_{2,2} + C_1 C_2, \qquad \qquad & b_{40} & = C_{4,4} + C_3 C_4,\notag
\end{alignat}
in addition to the following relations among the torsion and its derivatives:
\begin{alignat}{3}
B_1 & = -C_{1,\overline{0}}, \qquad \qquad & C_{1,3} & = C_1 C_4, \qquad \qquad & C_{3,1} & = C_2 C_3, \notag  \\
B_2 & = -C_{2,\overline{0}}, \qquad \qquad & C_{1,4} & = -C_1 C_3, \qquad \qquad & C_{3,2} & = -C_1 C_3, \label{B-relations} \\
B_3 & = -C_{3,\overline{0}}, \qquad \qquad & C_{2,3} & = C_2 C_4, \qquad \qquad & C_{4,1} & = C_2 C_4,  \notag \\
B_4 & = -C_{4,\overline{0}}, \qquad \qquad & C_{2,4} & = -C_2 C_3, \qquad \qquad & C_{4,2} & = -C_1 C_4. \notag
\end{alignat}

While we don't have an explicit coordinate representation for $\thetabar$, we can still explore the consequences of $d(d\thetabar) = 0$.  Computing $d(d\thetabar) \equiv 0$ modulo $\thetabar$ yields
\begin{gather*}
a_{11} = \frac{-2A_{2,1} + C_2 (A_2 + 2)}{2 A_2}, \qquad a_{42} = \frac{-2A_{2,2} - C_1 (A_2 + 2)}{2 A_2}, \\[0.2in]
A_{2,0} = A_2 (C_{3,3} + C_{4,4} - C_{1,1} - C_{2,2}).
\end{gather*}
Then computing $d(d\thetabar) \equiv 0$ modulo $\theta, \om^1, \om^2$ yields
\[ A_{2,\overline{0}} = C_{1,2} + C_{2,2} - A_2 (C_{3,3} + C_{4,4}) - \frac{C_1}{A_2} A_{2,1} - \frac{C_2}{A_2} A_{2,2} - C_3 A_{2,3} - C_4 A_{2,4}. \]
(Note that the fact that $\thetabar$ is a contact form implies that
$A_2$ cannot be zero.)

At this point, 
all coefficients in the structure equations \eqref{G0struct} have been expressed in terms of the functions $A_2, C_1, C_2, C_3, C_4$ and their first derivatives.  In addition, we have relations among the derivatives that amount to an overdetermined PDE system for these five functions on $\B$.  Necessary compatibility conditions for this system may be found by computing $d(dA_2) = d(dC_i) = 0$.  In particular, computing
\begin{align*}
d(dC_1) & \equiv d(dC_2) \equiv 0 \mod{\theta, \thetabar, \om^1, \om^2}, \\
d(dC_3) & \equiv d(dC_4) \equiv 0 \mod{\theta, \thetabar, \om^3, \om^4}
\end{align*}
yields
\begin{align}
C_{1,0} & = -A_2 C_{1,\overline{0}} - C_2 C_{1,2} + C_1 (-C_{1,1} + C_{3,3} + C_{4,4}), \notag \\
C_{2,0} & = -A_2 C_{2,\overline{0}} - C_1 C_{2,1} + C_2 (-C_{2,2} + C_{3,3} + C_{4,4}), \label{left-off-here} \\
C_{3,0} & = -C_{3,\overline{0}} - C_4 C_{3,4} + C_3 (C_{1,1} + C_{2,2} - C_{3,3}) , \notag \\
C_{4,0} & = -C_{4,\overline{0}} - C_3 C_{4,3} + C_4 (C_{1,1} + C_{2,2} - C_{4,4}). \notag
\end{align}

At this point, we have derived all the relations among the torsion functions on $\B$ and their derivatives that will be needed in order to prove that $(\Mbar, \Ibar)$ is Darboux-integrable after at most one prolongation.  The proof of Darboux-integrability is divided into two main cases.  In \S \ref{nondegenerate-case}, we prove Darboux-integrability under the assumption that the vectors $[C_1 \ \ C_2], \, [C_3 \ \ C_4]$ are both nonzero.  This case is further divided into three subcases, depending on the ranks of certain Frobenius systems that arise during the proof.  Precise statements of the results are contained in Propositions \ref{Case-1-prop}, \ref{Case-2-prop}, and \ref{Case-3-prop}.  In \S \ref{degenerate-case}, we prove Darboux-integrability under the assumption that exactly one of the vectors $[C_1 \ \ C_2], \, [C_3 \ \ C_4]$ vanishes; the precise result is contained in Proposition \ref{one-C-vector-zero-prop}.  As noted in \S \ref{wsetup}, it is not necessary to consider the case where $[C_1 \ \ C_2], \, [C_3 \ \ C_4]$ both vanish, since in that case both \MA/ systems are contact equivalent to the standard wave equation.

\subsection{Case 1: $[C_1 \ \ C_2], \, [C_3 \ \ C_4]$ are both nonzero}\label{nondegenerate-case}
Without loss of generality, we may assume that $C_2$ and $C_4$ are nonzero.  Consider the exterior derivatives of the ratios $\frac{C_1}{C_2}$ and $\frac{C_3}{C_4}$.  A straightforward computation shows that
\[
d\left(\frac{C_1}{C_2} \right) \equiv 0 \mod{\theta, \thetabar, \om^1, \om^2};
\]
therefore, $d\left(\frac{C_1}{C_2}\right)$ must lie in the last derived system  of $\calK_1 = \{\theta, \thetabar, \om^1, \om^2\}$---i.e., the largest integrable subsystem of $\calK_1$, denoted by $\calK_1^{(\infty)}$.  Similarly,
\[
d\left(\frac{C_3}{C_4} \right) \equiv 0 \mod{\theta, \thetabar, \om^3, \om^4},
\]
so $d\left(\frac{C_3}{C_4}\right)$ must lie in the last derived system of $\calK_2 = \{\theta, \thetabar, \om^3, \om^4\}$.

First, consider the system $\calK_1$.  In order to compute its first derived system $\calK_1^{(1)}$, we must find those 1-forms in $\calK_1$ whose exterior derivatives are zero modulo the linear span of the 1-forms in $\calK_1$.  To this end, we compute:
\[
\left. \begin{aligned}
d\theta & \equiv \om^3 \& \om^4 \\[0.1in]
d\thetabar & \equiv A_2\, \om^3 \& \om^4 \\[0.1in]
d\om^1 & \equiv C_1\, \om^3 \& \om^4 \\[0.1in]
d\om^2 & \equiv C_2\, \om^3 \& \om^4
\end{aligned} \right\} \mod{\calK_1}.
\]
Therefore,
$\calK_1^{(1)} = \{ \thetabar - A_2 \theta, \om^1 - C_1\theta, \om^2 - C_2\theta \}.$
Observe that
\[ \om^1 - C_1\theta = dX, \qquad  \om^2 - C_2\theta = dP. \]
Therefore, the rank 2 subsystem $\{\om^1 - C_1\theta, \om^2 - C_2\theta \} = \{dX, dP\}$ of $\calK_1^{(1)}$ is integrable, and the next derived system $\calK_1^{(2)}$ (i.e., the first derived system of $\calK_1^{(1)}$)  contains this rank 2 system.  The only question is whether, in fact, $\calK_1^{(2)} = \calK_1^{(1)}$---i.e., whether $\calK_1^{(1)}$ itself is integrable.  In either case, we will have $\calK_1^{(2)} = \calK_1^{(\infty)}$.  A computation shows that
\begin{equation}
 d(\thetabar - A_2 \theta) \equiv \theta \& [(A_{2,3} + A_2 C_4(A_2 - 1))\, \om^3 + (A_{2,4} - A_2 C_3(A_2 - 1))\, \om^4]  \mod \calK_1^{(1)}, \label{int-cond-1}
\end{equation}
so the rank of $\calK_1^{(\infty)}$ is either 3 or 2, depending on whether or not the 1-form in brackets vanishes.

Similarly, we can compute that
\[ \calK_2^{(1)} = \{ \thetabar - \theta, \om^3 - C_3\theta, \om^4 - C_4\theta \} = \{ \thetabar - \theta, dY, dQ \}. \]
So $\calK_2^{(\infty)}$ contains the rank 2 subsystem $\{\om^3 - C_3\theta, \om^4 - C_4\theta \} = \{dY, dQ\}$, and a computation shows that
\begin{equation}
 d(\thetabar - \theta) \equiv -\theta \& \left[\frac{(A_{2,1} + C_2(A_2 - 1))}{A_2}\, \om^1 + \frac{(A_{2,2} - C_1(A_2 - 1))}{A_2}\, \om^2 \right]  \mod \calK_2^{(1)}. \label{int-cond-2}
\end{equation}
So the rank of $\calK_2^{(\infty)}$ is either 3 or 2, 
depending on whether or not the 1-form in brackets vanishes.

Now we must divide into cases depending on the ranks of these derived systems.

\subsubsection{Case 1.1: $\calK_1^{(\infty)}$ and $\calK_2^{(\infty)}$ both have rank 3.}

In this case, we have the following result:

\begin{prop}\label{Case-1-prop}
If $[C_1 \ \ C_2], \, [C_3 \ \ C_4]$ are both nonzero and $\calK_1^{(\infty)}$ and $\calK_2^{(\infty)}$ both have rank 3, then the system $(\Mbar, \Ibar)$ is contact equivalent to the standard wave equation.
\end{prop}

\begin{proof}
By Theorem 2.1 of \cite{BGG}, it suffices to show that $(\Mbar, \Ibar)$ is Darboux-integrable, i.e., that each of the characteristic systems $\{\thetabar, \om^1, \om^2\}$ and $\{\thetabar, \om^3, \om^4\}$---which are well-defined on $\Mbar$ even though the 1-forms $\om^i$ are not---contains a rank 2 integrable subsystem.

The hypothesis that $\calK_1^{(\infty)}$ and $\calK_2^{(\infty)}$ both have rank 3 implies that the expressions \eqref{int-cond-1} and \eqref{int-cond-2} must both vanish identically; therefore,
$$
A_{2,1}  = -C_2(A_2-1), 
\qquad
A_{2,2}  = C_1(A_2-1), 
\qquad
A_{2,3}  = -A_2 C_4(A_2-1), 
\qquad
A_{2,4}  = A_2 C_3(A_2-1) .
$$
Using these conditions, a straightforward computation shows that
\begin{align*}
\{\thetabar, \om^1, \om^2\}^{(1)} & = \{ A_2 \om^1 - C_1 \thetabar, \, A_2 \om^2 - C_2 \thetabar \}, \\
\{\thetabar, \om^3, \om^4\}^{(1)} & = \{ A_2 \om^3 - C_3 \thetabar, \, A_2 \om^4 - C_4 \thetabar \},
\end{align*}
and that each of these derived systems is integrable.
\end{proof}

\subsubsection{Case 1.2: Exactly one of $\calK_1^{(\infty)}$ and $\calK_2^{(\infty)}$ has rank 3.}\label{middle-case}

Without loss of generality, we may assume that $\calK_2^{(\infty)}$ has rank 2 and is equal to $\{dY, dQ\}$, and that $\calK_1^{(\infty)}$ has rank 3.  It follows that \eqref{int-cond-1} vanishes identically and that \eqref{int-cond-2} does not.  Since all our results are local, we will assume that we are working on an open set where \eqref{int-cond-2} is nonzero.  The vanishing of \eqref{int-cond-1} implies that
\[ A_{2,3} = -A_2 C_4(A_2-1), \qquad
A_{2,4} = A_2 C_3 (A_2-1)  . \]

Recall that the function $\frac{C_3}{C_4}$ satisfies
\[ d \left(\frac{C_3}{C_4} \right) \in \calK_2^{(\infty)} = \{dY, dQ\}. \]
It follows that $\frac{C_3}{C_4}$ is a function of $Y$ and $Q$ alone.  Now consider the 1-form
\[ \widetilde{\om}^3 = \om^3 - \frac{C_3}{C_4}\, \om^4 = dY - \frac{C_3}{C_4} dQ. \]
This 1-form is contained in the span of $\om^3, \om^4$, and we have
\[ d\widetilde{\om}^3 \equiv 0 \mod \widetilde{\om}^3; \]
so $\widetilde{\om}^3$ is a multiple of an exact 1-form, say $\widetilde{\om}^3 = \lambda\, d\widetilde{Y}.$  Moreover, because $\widetilde{\om}^3$ is expressed solely in terms of $Y$ and $Q$, $\lambda$ and $\widetilde{Y}$ may be chosen to be functions depending only on $Y$ and $Q$, and which are therefore well-defined on $\scM$.  The crucial point here is that {\em there exists an exact 1-form in the span of $\{ \om^3, \om^4 \}$ which is well-defined on $\scM$.}  Then we have
\[ dY \& dQ = \widetilde{\om}^3 \& dQ = \lambda\, d\widetilde{Y} \& dQ = d\widetilde{Y} \& d\widetilde{Q}, \]
where
\[ \widetilde{Q}(\widetilde{Y}, Q) = \int_0^Q \lambda(\widetilde{Y},t)\, dt. \]
Since
\[ d\theta = dX \& dP + d\widetilde{Y} \& d\widetilde{Q} , \]
Pfaff's Theorem (see Ch. 1 of \cite{CFB}) implies that there exists a function $\widetilde{Z}$ on $\scM$ such that
\[ \theta = d\widetilde{Z} - P\, dX - \widetilde{Q}\, d\widetilde{Y}. \]

We can now repeat all our constructions starting with the coordinate system $(X,\widetilde{Y},\widetilde{Z},P,\widetilde{Q})$, but now our adapted coframing $\sigma$ will have
the additional property that $\om^3 = d\widetilde{Y}$ and $C_3 = 0$.  Thus we will drop the tildes and assume that $C_3 = 0$ for the remainder of this subsection.

\begin{prop}\label{Case-2-prop}
If $[C_1 \ \ C_2], \, [C_3 \ \ C_4]$ are both nonzero, $\calK_1^{(\infty)}$ has rank 3, and $\calK_2^{(\infty)}$ has rank 2, then the system $(\Mbar, \Ibar)$ is Monge-integrable, and it becomes Darboux-integrable after one prolongation.
\end{prop}

\begin{proof}
The same argument as that given in Case 1.1 shows that the characteristic system $\{\thetabar, \om^1, \om^2\}$ on $\Mbar$ contains a rank 2 integrable subsystem; therefore, $(\Mbar, \Ibar)$ is Monge-integrable.

In order to prove the second half of the Proposition, we will need to make use of relations among the second derivatives of the functions $A_2, C_1, C_2, C_4$.  These will be denoted as, e.g.,
\[ dA_{2,0} = A_{2,0 0}\theta + A_{2,0 \overline{0}}\thetabar + A_{2,0 1}\om^1 + A_{2,0 2}\om^2 + A_{2,0 3}\om^3 + A_{2,0 4}\om^4. \]
Note that, although (for example) the $A_{2,ij}$ are second derivatives of $A_2$, because we are working in a coframing rather than in coordinates, we cannot assume that $A_{2,ij} = A_{2,ji}.$

Computing $d(dA_2) \equiv 0 \mod{\{\thetabar, \om^3, \om^4\}}$ shows that
\begin{align*}
A_{2,10} & = A_2 (C_{4,41} - C_{1,11} - C_{2,21}  + C_2 (C_{1,1} + C_{2,2} - C_{4,4})) \\
& \qquad \qquad + A_{2,1} (C_{4,4} - C_{2,2} - C_1 C_2) + A_{2,2} (C_{2,1} - C_2^2),\\
A_{2,20} & = A_2 (C_{4,42} - C_{1,12} - C_{2,22} + C_1 (C_{4,4} - C_{1,1} - C_{2,2})) \\
& \qquad \qquad  + A_{2,1} (C_{1,2} + C_1^2) + A_{2,2} (C_{4,4} - C_{1,1} + C_1 C_2),   \\
A_{2,21} & =  A_{2,12} + (A_2-1)\left(C_{1,1} + C_{2,2}  - \frac{A_{2,1}}{A_2} C_1 - \frac{A_{2,2}}{A_2} C_2\right).
\end{align*}
Next, computing $d(dC_4) \equiv 0 \mod{\{\theta, \thetabar\}}$ shows that
\begin{align*}
C_{4,34} & = C_{4,43} + (A_2-1)C_{4,\overline{0}} + C_4(C_{1,1} + C_{2,2}), \\
C_{4,31} & = C_2 (C_{4,3} + C_4^2),\\
C_{4,32} & = -C_1 (C_{4,3} + C_4^2), \\
C_{4,41} & = C_2 C_{4,4}, \\
C_{4,42} & = -C_1 C_{4,4}.
\end{align*}
Now computing $d(dC_4) \equiv 0 \mod{\theta}$ shows that
\begin{align*}
C_{4, \overline{0} 1} & = C_4 C_{2,\overline{0}} + \frac{((A_2+1) C_1 - A_{2,1})}{A_2} C_{4,\overline{0}} , \\
C_{4, \overline{0} 2} & = -C_4 C_{1,\overline{0}} - \frac{((A_2+1) C_2 + A_{2,2})}{A_2} C_{4,\overline{0}} , \\
C_{4, 3\overline{0}} & = C_{4,\overline{0}3} - A_2 C_4 C_{4,\overline{0}}, \\
C_{4, 4\overline{0}} & = C_{4,\overline{0}4},
\end{align*}
and then computing $d(dC_4) \equiv 0 \mod{\thetabar}$ shows that
\begin{align*}
C_{4, 30} & = C_4 (C_{1,13} + C_{2,23} - C_{4,43} + C_{4,\overline{0}}) + (C_{4,3} - C_4^2)(C_{1,1} + C_{2,2}) - C_{4,\overline{0}3} , \\
C_{4, 40} & = C_4 (C_{1,14} + C_{2,24} - C_{4,44})  + C_{4,4}(C_{1,1} + C_{2,2}) - C_{4,\overline{0}4}, \\
C_{1,12} & = -C_{2,22} + (A_2-1) C_{1,\overline{0}} - C_1 C_{4,4} + \frac{((A_2 - 1) C_1 - A_{2,2})}{A_2 C_4} C_{4,\overline{0}},  \\
C_{2,21} & =   -C_{1,11} - (A_2-1) C_{2,\overline{0}} + C_2 C_{4,4} - \frac{((A_2 - 1) C_2 + A_{2,1})}{A_2 C_4} C_{4,\overline{0}} .
\end{align*}
Next, computing $d(dC_1) \equiv 0 \mod{\{\theta, \thetabar\}}$ yields
\begin{align*}
C_{1,13} & = C_4 (C_{1,1} + C_1 C_2), \\
C_{1,14} & = 0, \\
C_{1,23} & = C_4 (C_{1,2} - C_1^2), \\
C_{1,24} & = 0, \\
C_{1,21} & = -C_{2,22} - 2 C_1 C_{4,4} + 2 (A_2 - 1) C_{1,\overline{0}} + \frac{((A_2 - 1) C_1 - A_{2,2})}{A_2 C_4}   C_{4,\overline{0}}.
\end{align*}
Similarly, computing $d(dC_2) \equiv 0 \mod{\{\theta, \thetabar\}}$ yields
\begin{align*}
C_{2,13} & = C_4 (C_{42,1} + C_2^2), \\
C_{2,14} & = 0, \\
C_{2,23} & = C_4 (C_{2,2} - C_1 C_2), \\
C_{2,24} & = 0 ,\\
C_{2,12} & = -C_{1,11} + 2 C_2 C_{4,4} - 2 (A_2 - 1) C_{2,\overline{0}} - \frac{((A_2 - 1) C_2 + A_{2,1})}{A_2 C_4}   C_{4,\overline{0}}.
\end{align*}
Computing $d(dC_1) \equiv 0 \mod{\{\om^2, \om^3, \om^4\}}$ yields
\begin{align*}
C_{1,10} & = C_2 C_{2,22} - C_1 C_{1,11} - A_2 C_{1,\overline{0}1} + ((2 - A_2) C_2 - A_{2,1}) C_{1,\overline{0}} + C_{4,4} (C_{1,1} + 2 C_1 C_2 ) \\
& \qquad \qquad  - \frac{((A_2 - 1) C_1 - A_{2,2})}{A_2 C_{4}}  C_2 C_{4,\overline{0}},\\
C_{1,1\overline{0}} & = C_{1,\overline{0} 1} + \frac{(A_{2,1} - C_2)}{A_2} C_{1,\overline{0}},
\end{align*}
and computing $d(dC_2) \equiv 0 \mod{\{\om^1, \om^3, \om^4\}}$ yields
\begin{align*}
C_{2,20} & = C_1 C_{1,11} - C_2 C_{2,22} - A_2 C_{2,\overline{0}2} + ((A_2 - 2) C_1 - A_{2,2}) C_{2,\overline{0}} + C_{4,4} (C_{2,2} - 2 C_1 C_2 ) \\
& \qquad \qquad  + \frac{((A_2 - 1) C_2 + A_{2,1})}{A_2 C_{4}}  C_1 C_{4,\overline{0}},\\
C_{2,2\overline{0}} & = C_{2,\overline{0} 2} + \frac{(A_{2,2} + C_1)}{A_2} C_{2,\overline{0}},
\end{align*}
Now computing $d(dA_2) \equiv 0 \mod{\om^3}$ yields
\begin{align*}
A_{2,1\overline{0}} & = \frac{1}{A_2} \left[ -C_1 A_{2,11} - C_2 A_{2,12} + A_{2,1}(C_1 C_2 + C_{2,2}) + A_{2,2} (-C_{2,1} + C_2^2) \right]  \\
& \qquad \qquad  - (A_2 - 1) C_{2,\overline{0}} - C_2 (C_{1,1} + C_{2,2}) - ((A_2 - 2)C_2 + 2 A_{2,1}) C_{4,4}  \\
& \qquad \qquad - \frac{((A_2 - 1)C_2 + A_{2,1})}{A_2 C_4} C_{4,\overline{0}} ,\\
A_{2,2\overline{0}} & = \frac{1}{A_2} \left[ -C_1 A_{2,12} - C_2 A_{2,22} - A_{2,1}\left( C_{1,2} + \frac{C_1^2}{A_2} \right) + A_{2,2} \left( C_{1,1} - \frac{C_1 C_2}{A_2}\right) \right]  \\
& \qquad \qquad  + (A_2 - 1) C_{1,\overline{0}} + \frac{C_1}{A_2} (C_{1,1} + C_{2,2}) + ((A_2 - 2)C_1 - 2 A_{2,2}) C_{4,4}  \\
& \qquad \qquad + \frac{((A_2 - 1)C_1 - A_{2,2})}{A_2 C_4} C_{4,\overline{0}} ,
\end{align*}
and
\[ A_{2,14} = A_{2,24} = C_{4,44} = 0. \]

Finally, we need two additional relations which do not become apparent until we differentiate the equations for $dC_{4,3}$ and $dC_{4,4}$.  Computing $d(dC_{4,3}) \equiv 0 \mod{ \{\thetabar - \theta, \om^3, \om^4 - C_4 \theta\} }$ yields
\[ (C_4 C_{4,\overline{0}3} - ((A_2+1)C_4^2 + C_{4,3}) C_{4,\overline{0}}) \theta \& [(A_{2,1} + C_2(A_2 - 1))\, \om^1 + (A_{2,2} - C_1(A_2 - 1))\, \om^2]  . \]
Note that the right-hand factor is precisely \eqref{int-cond-2}, which we have assumed is nonzero.  Therefore,
\[ C_{4,\overline{0}3} = \frac{C_{4,\overline{0}}(C_{4,3} + (A_2+1) C_4^2)}{C_4} . \]
Precisely the same argument applied to $d(dC_{4,4})$ shows that
\[ C_{4,\overline{0}4} = \frac{C_{4,\overline{0}} C_{4,4}}{C_4}. \]

With these relations in hand, consider the characteristic system $\overline{\calK} = \{\thetabar, \om^3, \om^4\}$ of $\Ibar$---which is well-defined on $\Mbar$, even though $\om^4$ is not.  We need to show that after one prolongation, the corresponding characteristic system $\overline{\calK}'$ of the prolongation contains a rank 2 integrable subsystem.  In order to perform this computation, we need to construct a basis for $\overline{\calK}$ consisting of 1-forms which are well-defined on $\Mbar$.  Fortunately, $\thetabar$ and $\om^3$ are already well-defined on $\Mbar$.  For the remaining 1-form, it will be convenient to choose a 1-form which is contained in the first derived system $\overline{\calK}^{(1)} = \{\om^4 - C_4\, \thetabar, \om^3\}$.  To this end, introduce functions $\tau, g$ on $\B$ such that the 1-form
\[ \psi = e^{\tau}(\om^4 - C_4\, \thetabar - g\, \om^3) \]
is well-defined on $\Mbar$.  (The fact that $\overline{\calK}^{(1)}$ is well-defined on $\Mbar$ guarantees the existence of such functions.)
As before, we denote the derivatives of these functions by
\begin{align*}
d\tau & = \tau_0 \theta + \tau_{\overline{0}} \thetabar + \tau_1 \om^1 + \tau_2 \om^2 + \tau_3 \om^3 + \tau_4 \om^4 ,\\
dg & = g_0 \theta + g_{\overline{0}} \thetabar + g_1 \om^1 + g_2 \om^2 + g_3 \om^3 + g_4 \om^4,
\end{align*}
and similarly for second derivatives.

Because $\psi$ is well-defined on $\Mbar$,  $d\psi$ contains no terms involving $\theta$.  This, in turn, determines the partial derivatives $\tau_0, g_0$:
\begin{align*}
\tau_0 & = C_{4,4}, \\
g_0 & = C_4^2 - C_{4,3} - g C_{4,4}.
\end{align*}
We will also need to make use of relations among the second derivatives of $\tau, g$.  These are determined by computing $d(d\tau) = d(dg) = 0$; this is a straightforward computation, which we omit here for the sake of brevity.

We can define a {\em partial prolongation} $\Ibar'$ of $\Ibar$ on $\Mbar \times \R$ as follows.  (Note that Darboux-integrability of the partial prolongation implies Darboux-integrability of the full prolongation.)  Let $t$ be a new coordinate on the $\R$ factor; then the partial prolongation $\Ibar'$ is  generated by the 1-forms $\thetabar, \thetabar' = \psi - t \om^3,$ and the 2-form $\om^1 \& \om^2$.  Again, this system is well-defined on $\Mbar \times \R$, even though $\om^1, \om^2$ are not.

A straightforward computation shows that
 \[
d\thetabar' \equiv -\pi_1 \& \om^3
 \mod{\{\thetabar, \thetabar'\}},
\]
where
\[ \pi_1 = dt - \left(t \tau_1 - e^{\tau} g_1\right)\om^1 - \left(t \tau_2 - e^{\tau} g_2  \right)\om^2. \]
The corresponding characteristic system of $\Ibar'$ is
\[ \overline{\calK}' = \{\thetabar, \thetabar', \pi_1, \om^3 \}. \]

We will now compute the derived systems of $\overline{\calK}'$ and show that the second derived system $\overline{\calK}'^{(2)}$ is a Frobenius system of rank 2; this will complete the proof of the Proposition.  In order to compute the first derived system, we compute:
 \begin{equation*}
\left. \begin{aligned}
d\thetabar & \equiv \om^1 \& \om^2 \\[0.1in]
d\thetabar' & \equiv 0 \\[0.1in]
d\pi_1 & \equiv E \om^1 \& \om^2 \\[0.1in]
d\om^3 & \equiv 0
\end{aligned} \right\} \mod{\{\thetabar, \thetabar', \pi_1, \om^3\}}
\end{equation*}
(the last line following from $C_3=0$), where
\[
E = e^{\tau} C_4^2 + (e^{\tau}g_4 - t \tau_4) C_4 - e^{\tau} C_{4,3} - (e^{\tau} g + t) C_{4,4} + (e^{\tau} g_{\overline{0}} - t \tau_{\overline{0}}) .
\]
Let $\pi_2=\pi_1 - E \thetabar$, so that
\[ \overline{\calK}'^{(1)} = \{\thetabar', \pi_2, \om^3 \}. \]
Next we compute the derived system of $\overline{\calK}'^{(1)}$:
\[
\left. \begin{aligned}
d\thetabar' & \equiv \frac{C_4}{A_2} \thetabar \& \left[(A_{2,1} + C_2(A_2 - 1))\, \om^1 + (A_{2,2} - C_1(A_2 - 1))\, \om^2 \right]  \\[0.1in]
d\pi_2 & \equiv  \frac{F}{A_2} \thetabar \& \left[(A_{2,1} + C_2(A_2 - 1))\, \om^1 + (A_{2,2} - C_1(A_2 - 1))\, \om^2 \right] \\[0.1in]
d\om^3 & \equiv 0
\end{aligned} \right\} \mod{\{\thetabar', \pi_2, \om^3\}},
\]
where
\[  F = e^{\tau} C_4^2 + (e^{\tau}g_4 - t \tau_4) C_4 - e^{\tau} C_{4,3} - (e^{\tau} g + t) C_{4,4}. \]
Once again, we see the bracketed 1-form in \eqref{int-cond-2} appearing.  Since this 1-form is assumed to be nonzero, the derived system $\overline{\calK}'^{(2)}$ has rank 2 and is spanned by the forms $\om^3$ and
\[ \pi_3 =  C_4 \pi_2 - F\thetabar.  \]
Finally, another computation shows that
\[
\left. \begin{aligned}
d\pi_3 & \equiv  0 \\[0.1in]
d\om^3 & \equiv 0
\end{aligned} \right\} \mod{\{\pi_3, \om^3\}};
\]
therefore, $\overline{\calK}'^{(2)}$ is integrable.
\end{proof}

\subsubsection{Case 1.3: $\calK_1^{(\infty)}$ and $\calK_2^{(\infty)}$ both have rank 2.}
Now we assume that the bracketed 1-forms in both \eqref{int-cond-1} and \eqref{int-cond-2} are nonzero.  By the same argument as that given in the previous case, we may assume that $C_1 = C_3 = 0$, with $\om^1 = dX, \ \om^3 = dY$.

\begin{prop}\label{Case-3-prop}
If $[C_1 \ \ C_2], \, [C_3 \ \ C_4]$ are both nonzero and $\calK_1^{(\infty)}$ and $\calK_2^{(\infty)}$ both have rank 2, then the system $(\Mbar, \Ibar)$ becomes Darboux-integrable after one prolongation.
\end{prop}

\begin{proof}
The proof is very similar to that of Proposition \ref{Case-2-prop}.  We must now consider both characteristic systems
\[ \overline{\calK}_1 = \{\thetabar, \om^1, \om^2\}, \qquad \overline{\calK}_2 = \{\thetabar, \om^3, \om^4\} \]
of $\Ibar$.  As before, these systems are both well-defined on $\Mbar$, even though $\om^2$ and $\om^4$ are not.  We introduce functions $\rho, \tau, f, g$ on $\B$ such that the 1-forms
\[ \eta = e^{\rho}(A_2\om^2 - C_2\thetabar - f\om^1), \qquad \psi = e^{\tau}(\om^4 - C_4\thetabar - g\om^3) \]
are well-defined on $\Mbar$.  These forms have the property that
\[ \overline{\calK}_1^{(1)} = \{\eta, \om^1\}, \qquad \overline{\calK}_2^{(1)} = \{\psi, \om^3\}. \]

We construct the prolongation $\Ibar'$ of $\Ibar$ on $\Mbar \times \R^2$ as follows. Let $r, t$ be new coordinates on the $\R^2$ factor; then the prolongation $\Ibar'$ is  generated by the 1-forms $\thetabar, \thetabar_1 = \eta - r \om^1, \thetabar_2 = \psi - t \om^3$, and their exterior derivatives.

The remainder of the proof consists of applying the argument of Proposition \ref{Case-2-prop} to each of the characteristic systems $\overline{\calK}_1', \overline{\calK}_2'$ of the prolongation $\Ibar'$.  The argument varies only in the details of the calculations, and so we omit it for the sake of brevity.
\end{proof}

\subsection{Case 2: One of the $C$-vectors vanishes}\label{degenerate-case}

Without loss of generality, assume that $[C_1 \ \ C_2] = [0 \ \ 0]$, and that $C_4 \neq 0$.  By Proposition \ref{one-C-vector-zero}, it follows that $[B_1 \ \ B_2] = [0 \ \ 0]$ as well.

\begin{prop}\label{one-C-vector-zero-prop}
If $[C_1 \ \ C_2] = [0 \ \ 0]$, then the system $(\Mbar, \Ibar)$ is Monge-integrable, and it becomes Darboux-integrable after at most one prolongation. Furthermore, the \BT/ $\B \subset \scM \times \Mbar$ is holonomic.
\end{prop}

\begin{proof}
It follows from Proposition \ref{one-C-vector-zero} that $(\Mbar, \Ibar)$ is Monge-integrable; in fact, the characteristic system $\{\thetabar, \om^1, \om^2\}$ contains $\{\om^1, \om^2\} = \{dX, dP\}$ as a rank 2 integrable subsystem.

Now consider the other characteristic system $\calK = \{\thetabar, \om^3, \om^4\}$ of $\Ibar$.  One easily computes that the first derived system of $\calK$ is
\[ \calK^{(1)} = \{\om^3 - C_3 \thetabar, \om^4 - C_4 \thetabar\}. \]
In order to find the second derived system, we compute:
\[
\left.  \begin{aligned}
d(\om^3 - C_3 \thetabar) & \equiv \frac{C_3}{A_2} \thetabar \& (A_{2,1} \om^1 + A_{2,2}\om^2) \\
d(\om^4 - C_4 \thetabar) & \equiv \frac{C_4}{A_2} \thetabar \& (A_{2,1} \om^1 + A_{2,2}\om^2)
\end{aligned} \right\} \mod{\calK^{(1)}}.
\]
If $A_{2,1} = A_{2,2} = 0$, then $\calK^{(1)}$ is integrable; in this case, $\Ibar$ is Darboux-integrable and hence contact equivalent to the standard wave equation.  Therefore, we assume that $A_{2,1}$ and $A_{2,2}$ are not both zero.

In order to prove the second statement, we will construct a partial prolongation of $\Ibar$ and proceed as in \S \ref{middle-case}.  But first we need to derive relations among the derivatives of the torsion functions.

Picking up where we left off at \eqref{left-off-here}, consider $d(dC_3), d(dC_4)$.  Computing $d(dC_3) \equiv d(dC_4) \equiv 0 \mod{\{\theta - \thetabar, \om^3 - C_3 \thetabar, \om^4 - C_4 \thetabar\} }$ yields
\[ C_{3,\overline{0}} = C_{4,\overline{0}} = 0. \]
From \eqref{B-relations}, it follows that $B_3 = B_4 = 0$; therefore, the \BT/ is holonomic, as claimed.

We now have
\begin{align}
dC_{3} & = C_{3,3}(\om^3 - C_3 \theta) + C_{3,4} (\om^4 - C_4 \theta) = C_{3,3} dY + C_{3,4} dQ, \notag \\
dC_{4} & = C_{4,3}(\om^3 - C_3 \theta) + C_{4,4} (\om^4 - C_4 \theta) = C_{4,3} dY + C_{4,4} dQ. \label{dC4-exprn}
\end{align}
It follows that $C_3, C_4$ are functions of $Y$ and $Q$ alone.  Now the same argument as that given in \S \ref{middle-case} shows that we may assume $C_3 = 0$; moreover, $C_4$ remains a function of $Y$ and $Q$ alone when we do so.
Computing $d(dC_4) = 0$ yields the following relations among the second derivatives of $C_4$:
\[ C_{4,30} = - C_4 C_{4,34}, \qquad C_{4,40} = - C_4 C_{4,44}, \qquad C_{4,43} = C_{4,34}, \]
\[ C_{4,3\overline{0}} = C_{4,4\overline{0}} = C_{4,31} = C_{4,32} = C_{4,41} = C_{4,42} = 0. \]

Now consider the characteristic system $\calK = \{\thetabar, \om^3, \om^4\}$.  As we computed above (recalling that $C_3=0$), its first derived system is
\[ \calK^{(1)} = \{\om^3, \om^4 - C_4 \thetabar\}. \]
As in \S \ref{middle-case}, choose functions $g, \tau$ so that the 1-form
\[ \psi =  e^{\tau}(\om^4 - C_4\, \thetabar - g\, \om^3) \]
is well-defined on $\Mbar$, and construct the partial prolongation $\Ibar'$ of $\Ibar$ and the 1-form $\thetabar'$ as we did there.  Similar calculations to those of \S \ref{middle-case} show that the corresponding characteristic system $\overline{\calK}'$ of $\Ibar'$ has a rank 2 integrable subsystem.  This completes the proof.
\end{proof}

\section{Proof that Darboux implies B\"acklund}
\label{wprovebacklund}

\subsection{The non-Monge-integrable case}\label{NMIcase}
In this subsection $(\scM, \I)$ is assumed to be a hyperbolic \MA/ system which is
Darboux-integrable after one prolongation, but not Monge-integrable.
We will construct a canonical coframing associated to the prolongation.  We will then use this coframing to construct an integrable extension $\J$ of $(\scM, \I)$ in such a way that $\J$ defines a B\"acklund transformation between $(\scM, \I)$ and the standard wave equation $Z_{XY}=0$.

\begin{lemma}\label{NMIm5coframe} Near any point of $\scM$, there exists a coframing $(\theta,\pi_1,\pi_2,\eta^1,\eta^2)$
such that $\theta$ spans the 1-forms of $\I$, and the characteristic systems
$\calC_1, \calC_2$ of $\I$ have
derived flags
$$
\calC_1 = \{ \theta,\pi_1,\eta^1\} \supset \{\pi_1, \eta^1\} \supset \{\eta^1 \} =\calC_1^{(\infty)},
\qquad
\calC_2 = \{ \theta,\pi_2,\eta^2\} \supset \{\pi_2, \eta^2\} \supset \{\eta^2 \}=\calC_2^{(\infty)}.
$$
\end{lemma}
\begin{proof}
By a result of Jur\'a\v{s} \cite{Juras}, $(\scM,\I)$ is locally contact equivalent
to a system encoding a PDE of the form
\[ u_{xy}=F(x,y,u,p,q). \]
Thus, there are local coordinates
$x,y,u,p,q$ near the given point of $\scM$ such that $\I$ is generated by the 1-form
$\theta= du - p\, dx- q\, dy$ and
the 2-forms $(dp - F\,dy) \& dx$ and $(dq - F\,dx) \& dy$.  It is easy to verify that
the coframing given by $\theta$, $\eta^1 = dx$, $\eta^2 = dy$, $\pi_1 = dp - F\, dy - F_q \theta$,
and $\pi_2 = dq - F \,dx - F_p \theta$ has the properties claimed.
\end{proof}

In terms of the local coframing on $\scM$ given by the lemma, the
prolongation $(\scM', \I')$ is defined as follows:
let $\scM' = \scM \times \R^2$, with coordinates $r,t$ on the $\R^2$ factor,
and let $\I'$ be the Pfaffian system on $\scM'$
generated by $\theta$ and the forms
\begin{equation}\label{theta12def}
\theta_1 = \pi_1 - r \eta^1, \qquad \theta_2 = \pi_2 - t \eta^2.
\end{equation}

\begin{lemma}\label{n7coframe} Near any point of $\scM'$ there exists a coframing
$(\theta,\theta_1,\theta_2,\eta^1,\eta^2,\pi_3,\pi_4)$ such that
$\I'$ is generated by $\theta,\theta_1,\theta_2$, satisfying
\begin{equation}\label{nmi-Darboux-structure-0}
\begin{aligned}
d\theta & = -\theta_1 \& \eta^1 - \theta_2 \& \eta^2 \mod{\theta} 
\\
d\theta_1 & = -\pi_3 \& \eta^1 \mod{\theta, \theta_1} 
\\
d\theta_2 & = -\pi_4 \& \eta^2 \mod{\theta, \theta_2},
\end{aligned}
\end{equation}
with the derived flags of the characteristic systems of $\I'$ given by
\begin{gather*}
\Cprime1  = \{ \theta, \theta_1, \theta_2, \eta^1, \pi_3 \}
 \supset  \{ \theta, \theta_1, \eta^1, \pi_3 \}  \supset
 \{ \theta_1, \eta^1, \pi_3 \} \supset  \{ \eta^1, \pi_3 \}  = \Cprime1^{(\infty)},
\\
\Cprime2  = \{\theta, \theta_1, \theta_2, \eta^2, \pi_4\}
\supset \{ \theta, \theta_2, \eta^2, \pi_4 \} \supset
\{\theta_2, \eta^2, \pi_4\} \supset  \{\eta^2, \pi_4\} = \Cprime2^{(\infty)}.
\end{gather*}
\end{lemma}

\begin{proof}Let $\theta,\eta^1=dx, \eta^2=dy$ be part of the local coframing
on $\scM$ (pulled back to $\scM')$ constructed in the proof of Lemma \ref{NMIm5coframe},
and let $\theta_1,\theta_2$ be defined as in \eqref{theta12def}. Then
\begin{align*}
d\theta_1 &\equiv -(dr - (D_x F) dy) \& dx \quad \mod{\theta,\theta_1}, \\
d\theta_2 &\equiv -(dt - (D_y F) dx) \& dy \quad \mod{\theta,\theta_2},
\end{align*}
where
\begin{align*}
D_x F &= F_x + F_u p + F_p r + F_q F,\\
D_y F &= F_y + F_u q + F_p F + F_q t.
\end{align*}
For the moment, let $\pi_3 = dr - (D_x F)dy$.
Because $d\pi_3  \equiv 0$ modulo $dx, \pi_3, \theta,\theta_1,\theta_2$,
it follows that $\pi_3$ lies in $\Cprime1^{(1)}$.  Moreover, we may subtract a multiple of
$\theta$ from $\pi_3$ to ensure that $\pi_3$ lies in $\Cprime1^{(2)}$.

Next, we prove that the last derived system of $\Cprime1$ has rank 2, rather than rank 3.
(A similar argument applies to $\Cprime2$.)  Suppose that $\Cprime1^{(2)} = \{\theta_1, \eta^1, \pi_3\}$ is integrable.  From \eqref{theta12def}, it is clear that this is equivalent to the statement that $\{\pi_1, \eta^1, \pi_3\}$ is integrable---i.e., that $d\pi_1 \equiv d\eta^1 \equiv d\pi_3 \equiv 0$ modulo $\pi_1, \eta^1, \pi_3$.  But $\pi_1$ and $\eta^1$ are both well-defined on $\scM$, so their exterior derivatives do not involve $\pi_3$.  It follows that $d\pi_1 \equiv d\eta^1 \equiv 0$ modulo $\pi_1, \eta^1$, and $\calC_1^{(1)} = \{\pi_1, \eta^1\}$ is integrable, contrary to the hypothesis that $(\scM, \I)$ is not Monge-integrable.


\end{proof}

The conditions in Lemma \ref{n7coframe} are preserved by changes of coframing
of the form
\begin{equation} \label{Darboux-freedom}
\begin{bmatrix}
\tilde{\theta}\\[0.1in] \tilde{\theta}_1 \\[0.1in] \tilde{\theta}_2 \\[0.1in] \tilde{\eta}^1 \\[0.1in] \tilde{\eta}^2 \\[0.1in] \tilde{\pi}_3 \\[0.1in] \tilde{\pi}_4
\end{bmatrix} =
\begin{bmatrix}
c & 0 & 0 & 0 & 0 & 0 & 0 \\[0.1in]
0 & a_1 c & 0 & 0 & 0 & 0 & 0  \\[0.1in]
0 & 0 & a_2 c & 0 & 0 & 0 & 0  \\[0.1in]
0 & 0 & 0 & a_1^{-1} & 0 & 0 & 0  \\[0.1in]
0 & 0 & 0 & 0 & a_2^{-1} & 0 & 0   \\[0.1in]
0 & 0 & 0 & b_1 & 0 & a_1^2 c & 0  \\[0.1in]
0 & 0 & 0 & 0 & b_2 & 0 & a_2^2 c
\end{bmatrix}^{-1}
\begin{bmatrix}
\theta\\[0.1in] \theta_1 \\[0.1in] \theta_2 \\[0.1in] \eta^1 \\[0.1in] \eta^2 \\[0.1in] \pi_3 \\[0.1in] \pi_4
\end{bmatrix},
\end{equation}
with $a_1, a_2, c \neq 0$.  Let $G\subset GL(7,\R)$ be the group
of such transformations, and let $\calP$ be the $G$-structure on $\scM'$
of which the coframing of Lemma \ref{n7coframe} is a section.

After absorbing as much torsion as possible and differentiating to uncover relations among the torsion,
$\calP$ has structure equations
\begin{multline}\hskip -.2in
\begin{bmatrix} d\theta\\[0.1in] d\theta_1\\[0.1in] d\theta_2 \\[0.1in] d\eta^1\\[0.1in] d\eta^2 \\[0.1in] d\pi_3 \\[0.1in] d\pi_4
\end{bmatrix} =-
{\setlength\arraycolsep{0.5pt}
\left[  \begin{array}{c c c c c c c}
\gamma & 0 & 0 & 0 & 0 & 0 & 0  \\[0.1in]
0 & \ \gamma + \alp_1 & 0 & 0 & 0 & 0 & 0 \\[0.1in]
0 & 0 & \gamma + \alp_2 & 0 & 0 & 0 & 0 \\[0.1in]
0 & 0 & 0 & -\alp_1 & 0 & 0 & 0 \\[0.1in]
0 & 0 & 0 & 0 & -\alp_2 & 0 & 0 \\[0.1in]
0 & 0 & 0 & \bet_1 & 0 & \gamma + 2\alp_1 & 0 \\[0.1in]
0 & 0 & 0 & 0 & \bet_2 & 0 & \gamma + 2\alp_2
\end{array}\right]
}
 \&
\begin{bmatrix} \theta\\[0.1in] \theta_1\\[0.1in] \theta_2 \\[0.1in] \eta^1\\[0.1in] \eta^2 \\[0.1in] \pi_3 \\[0.1in] \pi_4
\end{bmatrix}
- \begin{bmatrix}
\theta_1 \& \eta^1 \! + \! \theta_2 \& \eta^2 \\[0.1in]
\pi_3 \& \eta^1 \! + \! (A_2 \theta_2 \! + \! B_2 \eta^2) \& \theta \\[0.1in]
\pi_4 \& \eta^2 \! + \! (A_1 \theta_1 \! + \! B_1 \eta^1) \& \theta \\[0.1in]
0 \\[0.1in]
0 \\[0.1in]
2C_1 \theta_1 \& \pi_3   \\[0.1in]
2C_2 \theta_2 \& \pi_4
\end{bmatrix}\hskip -2pt
.
\end{multline}

Because of the dimensions of the derived flags of the characteristic systems (given
in Lemma \ref{n7coframe}),  $A_1, B_1$ are not both zero, and $A_2, B_2$ are not both zero.
Furthermore, we can choose a local section $\sigma: \scM' \to \calP$ satisfying the conditions that $\eta^1 = dx$, $\eta^2 = dy$, and the forms $\pi_3, \pi_4$ are integrable; i.e.,
\[ d\pi_3 \equiv 0 \mod{\pi_3}, \qquad d\pi_4 \equiv 0 \mod{\pi_4}. \]
To see why, note that $\{ \eta^1, \pi_3 \}$ is a Frobenius system, and so it is spanned locally by two exact 1-forms.  Thus we can adjust $\pi_3$ by adding multiples of $\eta^1$ in order to make it a multiple of an exact form.   Similarly, we can add multiples of $\eta^2$ to  $\pi_4$ in order to make $\pi_4$ a multiple of an exact form.  However, we cannot independently scale $\pi_3$ and $\pi_4$ to make both of them exact.

This choice of section is not unique; it is determined up to a transformation of the form
\begin{equation}
\begin{bmatrix}
\tilde{\theta}\\[0.1in] \tilde{\theta}_1 \\[0.1in] \tilde{\theta}_2 \\[0.1in] \tilde{\eta}^1 \\[0.1in] \tilde{\eta}^2 \\[0.1in] \tilde{\pi}_3 \\[0.1in] \tilde{\pi}_4
\end{bmatrix} =
\begin{bmatrix}
c & 0 & 0 & 0 & 0 & 0 & 0 \\[0.1in]
0 & c & 0 & 0 & 0 & 0 & 0  \\[0.1in]
0 & 0 & c & 0 & 0 & 0 & 0  \\[0.1in]
0 & 0 & 0 & 1 & 0 & 0 & 0  \\[0.1in]
0 & 0 & 0 & 0 & 1 & 0 & 0   \\[0.1in]
0 & 0 & 0 & 0 & 0 &  c & 0  \\[0.1in]
0 & 0 & 0 & 0 & 0 & 0 &  c
\end{bmatrix}^{-1}
\begin{bmatrix}
\theta\\[0.1in] \theta_1 \\[0.1in] \theta_2 \\[0.1in] \eta^1 \\[0.1in] \eta^2 \\[0.1in] \pi_3 \\[0.1in] \pi_4
\end{bmatrix}
\end{equation}
with $c \neq 0$.  However, we can make the choice of $\sigma$ unique (albeit slightly non-canonical) as follows: since $\pi_3$, $\pi_4$ are integrable 1-forms, we must have
\[ \pi_3 = e^g d\xi_1, \qquad \pi_4 = e^h d\xi_2 \]
for some functions $\xi_1, \xi_2, f, g$ on $\scM'$.  Using the remaining scaling freedom, we can arrange that $h=-g$; the resulting coframing $\sigma: \calP \to \scM$ is uniquely determined.

When we pull back the structure equations via $\sigma$, the pseudoconnection forms $\alp_1$, $\alp_2$, $\beta_1$, $\beta_2$, $\gamma$ become semi-basic.  By making use of the remaining ambiguity in these forms and the conditions imposed thus far on the coframing, we can assume that
\begin{align*}
\alp_1 & = (D_1 + E_1)\eta^1 \\
\alp_2 & = (D_2 + E_2)\eta^2 \\
\beta_1 & = 2 D_1 \pi_3 \\
\beta_2 & = 2 D_2 \pi_4 \\
\gamma & = -C_1 \theta_1 - C_2 \theta_2 - E_1 \eta^1 - E_2 \eta^2 + F_1 \pi_3 + F_2 \pi_4
\end{align*}
for some functions $C_i, D_i, E_i, F_i$.
Then the structure equations for this coframing become:
\begin{align}
d\theta & = \theta \& (-C_1 \theta_1 - C_2 \theta_2 - E_1 \eta^1 - E_2 \eta^2 + F_1 \pi_3 + F_2 \pi_4)  - \theta_1 \& \eta^1 - \theta_2 \& \eta^2   \notag \\
d\theta_1 & = \theta_1 \& (-C_2 \theta_2 - E_2 \eta^2 + F_1 \pi_3 + F_2 \pi_4)  + D_1 \theta_1 \& \eta^1  - \pi_3 \& \eta^1 + \theta \& (A_2 \theta_2 + B_2 \eta^2)   \notag \\
d\theta_2 & = \theta_2 \& (-C_1 \theta_1 - E_1 \eta^1 + F_1 \pi_3 + F_2 \pi_4)  + D_2 \theta_2 \& \eta^2  - \pi_4 \& \eta^2 + \theta \& (A_1 \theta_1 + B_1 \eta^1) \notag \\
d\eta^1 & = 0   \label{Darboux-section-structure} \\
d\eta^2 & =  0  \notag \\
d\pi_3 & = \pi_3 \& (C_1 \theta_1 - C_2 \theta_2 + E_1 \eta^1 - E_2 \eta^2 + F_2 \pi_4)    \notag \\
d\pi_4 & = \pi_4 \& (-C_1 \theta_1 + C_2 \theta_2 - E_1 \eta^1 + E_2 \eta^2 + F_1 \pi_3).   \notag
\end{align}
(Note that these torsion functions are completely unrelated to those in \S \ref{wsetup} and \S \ref{wprovedarboux}.)

As in \S \ref{wprovedarboux}, we will need to compute relations among the derivatives of the torsion functions in order to show that $(\scM, \I)$ has a B\"acklund transformation to the wave equation.  We begin by differentiating the structure equations \eqref{Darboux-section-structure}.  Using notation similar to that in \S \ref{wprovedarboux}, we denote derivatives as, e.g.,
\[ dA_1 = A_{1,0} \theta + A_{1,1} \theta_1 + A_{1,2} \theta_2 + A_{1,3} \eta^1 + A_{1,4} \eta^2 + A_{1,5} \pi_3 + A_{1,6} \pi_4. \]
(Note that since this coframing is defined on a different manifold from that in \S \ref{wprovedarboux}, the indexing of the derivatives is different as well.)

Computing $d(d\theta) = d(d\theta_1) = d(d\theta_2) = d(d\pi_3)  = d(d\pi_4) = 0$ yields the following equations for the derivatives of the torsion functions:
\begin{align}
dA_1 & = A_{1,0} \theta + A_{1,1} \theta_1 -3 A_1 C_2 \theta_2 + A_{1,3} \eta^1 - A_1(D_2 + 2 E_2) \eta^2 + A_1 F_3 \pi_3 + A_1 F_2 \pi_4 \notag \\
dA_2 & = A_{2,0} \theta - 3 A_2 C_1 \theta_1 + A_{2,2} \theta_2 - A_2(D_1 + 2 E_1) \eta^1 + A_{2,4} \eta^2 + A_2 F_3 \pi_3 + A_2 F_2 \pi_4 \notag \\
dB_1 & = B_{1,0} \theta + (A_{1,3} - A_1 D_1) \theta_1 - 2 B_1 C_2 \theta_2 + B_{1,3} \eta^1 - B_1(D_2 + E_2) \eta^2 + A_1 \pi_3 \notag \\
dB_2 & = B_{2,0} \theta - 2 B_2 C_1 \theta_1 + (A_{2,4} - A_2 D_2) \theta_2  - B_2(D_1 + E_1) \eta^1 + B_{2,4} \eta^2 + A_2 \pi_4 \notag \\
dC_1 & = A_1 C_2 \theta + C_{1,1} \theta_1 - C_1 C_2 \theta_2 + C_{1,3} \eta^1 - (\tfrac{1}{2} A_1 + C_1 E_2) \eta^2 + C_{1,5} \pi_3 + C_1 F_2 \pi_4 \notag \\
dC_2 & = A_2 C_1 \theta - C_1 C_2 \theta_1 + C_{2,2} \theta_2 - (\tfrac{1}{2} A_2 + C_2 E_1) \eta^1 + C_{2,4} \eta^2 + C_2 F_1 \pi_3 + C_{2,6} \pi_4 \label{dtorsion} \\
dD_1 & = B_1 C_2 \theta + D_{1,1} \theta_1 + \tfrac{3}{2} A_2 \theta_2 + D_{1,3} \eta^1 + \tfrac{1}{2}(3 B_2 - B_1) \eta^2 + (2 C_1 - E_{1,5}) \pi_3 \notag \\
dD_2 & = B_2 C_1 \theta + \tfrac{3}{2} A_1 \theta_1 + D_{2,2} \theta_2 + \tfrac{1}{2}(3 B_1 - B_2) \eta^1 + D_{2,4} \eta^2 + (2 C_2 - E_{2,6}) \pi_4 \notag \\
dE_1 & = B_1 C_2 \theta + (C_{1,3} - C_1 D_1) \theta_1 + \tfrac{1}{2} A_2 \theta_2 + E_{1,3} \eta^1 + \tfrac{1}{2}(B_2 - B_1) \eta^2 + E_{1,5} \pi_3 \notag \\
dE_2 & = B_2 C_1 \theta + \tfrac{1}{2} A_1 \theta_1 + (C_{2,4} - C_2 D_2) \theta_2 + \tfrac{1}{2}(B_1 - B_2) \eta^1 + E_{2,4} \eta^2 + E_{2,6} \pi_4 \notag \\
dF_1 & = (2 C_1 F_1 - C_{1,5}) \theta_1 - C_2 F_1 \theta_2 + (C_1 + E_1 F_1 - E_{1,5}) \eta^1 - E_2 F_1 \eta^2 + F_{1,5} \pi_3 + F_{1,6} \pi_4\notag \\
dF_2 & = -C_1 F_2 \theta_1 + (2 C_2 F_2 - C_{2,6}) \theta_2 - E_1 F_2 \eta^1 + (C_2 + E_2 F_2 - E_{2,6}) \eta^2 + F_{1,6} \pi_3 + F_{2,6} \pi_4. \notag
\end{align}

Because $A_1$ appears as a derivative of $B_1$, and $A_1, B_1$ cannot vanish simultaneously, $B_1$ cannot vanish on any open set in $\scM'$.  In fact, $B_1$ cannot vanish identically on any fiber of the projection $\scM' \to \scM$, and the same is true of $B_2$.   Henceforth we restrict to the dense open set in $\scM'$ where $B_1, B_2$ are both nonzero, and note that this set surjects onto $\scM$.

We may obtain further relations among the derivatives of the torsion functions by differentiating equations \eqref{dtorsion}.  Computing $d(dA_1) \equiv d(dB_1) \equiv 0$ modulo $\theta, \theta_1, \eta^1, \pi_4$ yields
\[
A_{1,0}  = A_1(C_{2,4} - D_{2,2} - C_2 D_2), \qquad
B_{1,0}  = B_1(C_{2,4} - D_{2,2} - C_2 D_2).
\]
Then computing $d(dB_1) \equiv 0$ modulo $\theta, \theta_1, \eta^1$ yields
$ C_{2,6} = C_2 F_2. $
Similar considerations of $d(dA_2)$ and $d(dB_2)$ show that
\begin{gather*}
A_{2,0}  = A_2(C_{1,3} - D_{1,1} - C_1 D_1) , \qquad
B_{2,0}  = B_2(C_{1,3} - D_{1,1} - C_1 D_1), \qquad
C_{1,5} = C_1 F_1.
\end{gather*}

It will now be convenient to derive several equations and solve them simultaneously.  First, $d(dC_1) \equiv 0$ modulo $\theta_1, \eta^1$ implies that
\begin{gather}
A_1(C_{2,2} - C_2^2) = A_2(C_{1,1} - C_1^2) \label{eqn1forH} \\
A_1(3C_{2,4} - D_{2,2} - 3 C_2 D_2) = 2 B_2 (C_{1,1} - C_1^2).
\end{gather}
Additionally, $d(dC_2) \equiv 0$ modulo $\theta_2, \eta^2$ implies that
\begin{equation}
A_2(3C_{1,3} - D_{1,1} - 3 C_1 D_1) = 2 B_1 (C_{2,2} - C_2^2).
\end{equation}
Finally, $d(dD_1) \equiv 0$ modulo $\theta_1, \eta^1, \pi_3$ implies that
\begin{equation}
B_1(3C_{2,4} - D_{2,2} - 3 C_2 D_2) = B_2(3C_{1,3} - D_{1,1} - 3 C_1 D_1). \label{eqn4forH}
\end{equation}
The general solution to equations \eqref{eqn1forH}-\eqref{eqn4forH} is most easily expressed in terms of a new torsion function $H$, such that
\begin{align*}
C_{1,1} &= C_1^2 - A_1 H, &
D_{1,1} &= 3(C_{1,3} - C_1 D_1) + 2 B_1 H, \\
C_{2,2} &= C_2^2 - A_2 H, &
D_{2,2} &= 3(C_{2,4} - C_2 D_2) + 2 B_2 H.
\end{align*}

Next, we need information about the derivatives of $C_{1,3}$ and $C_{2,4}$.  Computing $d(dC_1) \equiv 0$ modulo $\theta_1$ yields
\begin{align*}
dC_{1,3} & = (C_2 A_{1,3} - \tfrac{1}{2} A_1 A_2 + B_1 C_1 C_2) \theta + C_{1,31} \theta_1 + (\tfrac{1}{2} A_2 C_1 - C_2 C_{1,3}) \theta_2 \\
& \qquad  + C_{1,33} \eta^1  - \tfrac{1}{2}(A_{1,3} + 2 E_2 C_{1,3} + C_1(B_1 - B_2)) \eta^2 \\
& \qquad \qquad + (F_1 C_{1,3} - C_1 E_{1,5} + 2 C_1^2 - A_1 H) \pi_3  + F_2 C_{1,3} \pi_4,
\end{align*}
and computing $d(dC_2) \equiv 0$ modulo $\theta_2$ yields
\begin{align*}
dC_{2,4} & = (C_1 A_{2,4} - \tfrac{1}{2} A_1 A_2 + B_2 C_1 C_2) \theta + (\tfrac{1}{2} A_1 C_2 - C_1 C_{2,4}) \theta_1 + C_{2,42} \theta_2 \\
& \qquad  - \tfrac{1}{2}(A_{2,4} + 2 E_1 C_{2,4} + C_2(B_2 - B_1)) \eta^1
 + C_{2,44} \eta^2 \\
 & \qquad \qquad  + F_1 C_{2,4} \pi_3 + (F_2 C_{2,4} - C_2 E_{2,4} + 2 C_2^2 - A_2 H) \pi_4 .
\end{align*}

Now, computing $d(dC_1) = d(dC_2) = 0$, $d(dD_1) \equiv 0$ modulo $\eta^1, \pi_3$, and $d(dD_2) \equiv 0$ modulo $\eta^2, \pi_4$ yields four different expressions for $dH$.  Taking linear combinations of these expressions shows that
\begin{gather}
A_1 (A_{1,3} - A_1 D_1) = B_1 (A_{1,1} - A_1 C_1) \label{J1eq} \\
A_2 (A_{2,4} - A_2 D_2) = B_2 (A_{2,2} - A_2 C_2) \label{J2eq}.
\end{gather}
Equations \eqref{J1eq} and \eqref{J2eq} may be solved by introducing new torsion functions $J_1, J_2$,
such that
\begin{align*}
A_{1,1} & = A_1 C_1 + A_1 J_1 &
A_{2,2} & = A_2 C_2 + A_2 J_2 \\
A_{1,3} & = A_1 D_1 + B_1 J_1 &
A_{2,4} & = A_2 D_2 + B_2 J_2.
\end{align*}
Then the various expressions for $dH$ may be combined to show that
\begin{gather*}
C_{1,31} = 2 C_1 C_{1,3} - C_1^2 D_1 + A_1 C_2 - \tfrac{1}{2} A_1 J_2 - (A_1 D_1 + B_1 J_1) H \\
C_{2,42} = 2 C_2 C_{2,4} - C_2^2 D_2 + A_2 C_1 - \tfrac{1}{2} A_2 J_1 - (A_2 D_2 + B_2 J_2) H,
\end{gather*}
and
\begin{multline*}
dH = \left(H(C_{1,3} - C_1 D_1 + B_1 H + C_{2,4} - C_2 D_2 + B_2 H) - \tfrac{1}{2} C_1 J_2 + C_2 J_1\right) \theta \\ + C_1 H \theta_1 + C_2 H \theta_2 + (D_1 H + \tfrac{1}{2} J_2) \eta^1 + (D_2 H + \tfrac{1}{2} J_1) \eta^2
+ F_1 H \pi_3 + F_2 H \pi_4.
\end{multline*}
They also imply the relation
\begin{equation}
H(C_{1,3} - C_1 D_1 + B_1 H) - \tfrac{1}{2} C_1 J_2 = H(C_{2,4} - C_2 D_2 + B_2 H) - \tfrac{1}{2} C_2 J_1. \label{zero1}
\end{equation}
The equations for $dA_1, dA_2$ now take the form:
\begin{multline}
dA_1 = -2A_1(C_{2,4} - C_2 D_2 + B_2 H)  \theta + A_1(C_1 + J_1) \theta_1 -3 A_1 C_2 \theta_2 \\
+ (A_1 D_1 + B_1 J_1) \eta^1  - A_1(D_2 + 2 E_2) \eta^2 + A_1 F_3 \pi_3 + A_1 F_2 \pi_4 \label{betterdA1eq}
\end{multline}
\begin{multline}
dA_2 =  -2A_2(C_{1,3} - C_1 D_1 + B_1 H) \theta - 3 A_2 C_1 \theta_1 + A_2(C_2 + J_2) \theta_2 \\
 - A_2(D_1 + 2 E_1) \eta^1 + (A_2 D_2 + B_2 J_2) \eta^2 + A_2 F_3 \pi_3 + A_2 F_2 \pi_4. \label{betterdA2eq}
\end{multline}

Now, computing $d(dB_1) \equiv 0$ modulo $\theta_1, \eta^1$ yields
\[ C_{2,44} = (D_2 - E_2)C_{2,4} + C_2 D_{2,4} - H B_{2,4} - (D_2 + E_2) B_2 H + C_2 D_2 E_2 + B_2 (C_1 - J_1). \]
Similarly, computing $d(dB_2) \equiv 0$ modulo $\theta_2, \eta^2$ yields
\[ C_{1,33} = (D_1 - E_1)C_{1,3} + C_1 D_{1,3} - H B_{1,3} - (D_1 + E_1) B_1 H + C_1 D_1 E_1 + B_1 (C_2 - J_2). \]
Next, computing $d(dA_1) = 0$ shows that
\begin{align*}
dJ_1 & = 4 A_1 C_2 \theta + (C_1 J_1 + A_1 K_1) \theta_1 - C_2 J_1 \theta_2 \\
& \qquad + (2C_{1,3} - 2C_1 D_1 + 2B_1 H + 2C_{2,4} - 2C_2 D_2 + 2B_2 H + D_1 J_1 + B_1 K_1) \eta^1 \\
& \qquad \qquad  - (2 A_1 + E_2 J_1) \eta^2 + F_1 J_1 \pi_3 + F_2 J_1 \pi_4
\end{align*}
for some function $K_1$.  Similarly, computing $d(dA_2) = 0$ shows that
\begin{align*}
dJ_2 & = 4 A_2 C_1 \theta - C_1 J_2 \theta_1 + (C_2 J_2 + A_2 K_2) \theta_2 - (2 A_2 + E_1 J_2) \eta^1\\
& \qquad + (2C_{1,3} - 2C_1 D_1 + 2B_1 H + 2C_{2,4} - 2C_2 D_2 + 2B_2 H + D_2 J_2 + B_2 K_2) \eta^2 \\
& \qquad \qquad   + F_1 J_2 \pi_3 + F_2 J_2 \pi_4
\end{align*}
for some function $K_2$.  But now computing $d(dH) \equiv 0$ modulo $\theta, \theta_1, \theta_2$ yields
\[ B_1 (K_1 + 4 H) = B_2 (K_2 + 4 H). \]
It follows that
\begin{align*}
K_1 & = -4H + B_2 M \\
K_2 & = -4H + B_1 M
\end{align*}
for some function $M$.  Computing $d(dH) \equiv 0$ modulo $\theta_2, \eta^1$ and $d(dH) \equiv 0$ modulo $\theta_1, \eta^2$ shows that
\begin{equation}
 A_1 M = A_2 M = 0. \label{Meqn}
\end{equation}

\begin{claim*}
$M = 0$.
\end{claim*}

\begin{pf}
Suppose not. Then by \eqref{Meqn}, $A_1 = A_2 = 0.$  Therefore, equations \eqref{betterdA1eq} and \eqref{betterdA2eq} reduce to
\begin{gather*}
 0 = dA_1 = B_1 J_1 \eta^1 \\
 0 = dA_2 = B_2 J_2 \eta^2.
\end{gather*}
Since $B_1, B_2$ are nonzero, it follows that $J_1 = J_2 = 0$.  Then
\begin{gather*}
0 = dJ_1 = (2C_{1,3} - 2C_1 D_1 - 2B_1 H + 2C_{2,4} - 2C_2 D_2 + 2B_2 H + B_1 B_2 M) \eta^1 \\
0 = dJ_2 = (2C_{1,3} - 2C_1 D_1 + 2B_1 H + 2C_{2,4} - 2C_2 D_2 - 2B_2 H + B_1 B_2 M) \eta^2.
\end{gather*}
Subtracting the two coefficients above yields
\[ 4 (B_2 - B_1) H = 0, \]
so either $B_1 = B_2$ or $H=0$.

First suppose that $B_1= B_2$. Then
\[ 0 = d(B_2 - B_1) \equiv -2B_1 (C_1 \theta_1 - C_2 \theta_2) \mod{\theta, \eta^1, \eta^2}, \]
so $C_1 = C_2 = 0.$  It follows that $C_{1,3} = C_{2,4} = 0$ as well.  But now
\[ 0 = dJ_1 = B_1^2 M, \]
so $M=0$, as desired.

Now suppose that $H=0$.  Computing $d(dC_{1,3}) - C_1 d(dD_1) \equiv 0$ modulo $\theta$ yields
\[ B_1^2 B_2 M = 0, \]
so $M=0$ in this case as well.
\end{pf}

Finally, computing $d(dH) = 0$, keeping the relation \eqref{zero1} in mind, yields two additional relations:
\begin{gather}
2C_2 (2C_{1,3} - 2 C_1 D_1 + 2 B_1 H) + (2C_2 - J_2) (2C_{2,4} - 2 C_2 D_2 + 2 B_2 H) + A_2 (4 C_1 - J_1) = 0 \label{zero2} \\
(2C_1 - J_1)(2C_{1,3} - 2 C_1 D_1 + 2 B_1 H) + 2 C_1 (2C_{2,4} - 2 C_2 D_2 + 2 B_2 H) + A_1 (4 C_2 - J_2) = 0. \label{zero3}
\end{gather}
We now have all the relations that will be needed for the involutivity calculation below.

\bigskip

Now suppose that $\Mbar=\R^5$ carries a \MA/ system $\Ibar$ representing
the wave equation $Z_{XY}=0$, generated algebraically by the contact form
\begin{equation}\label{defthetabar}
\thetabar = dZ - P\, dX - Q\, dY
\end{equation}
and the 2-forms $dP \& dX$ and $dQ \& dY$.
If there were a B\"acklund transformation $\B \subset \scM \times \Mbar$,
then $Z$ would be a local coordinate on the fibers of $\B \to \scM$
and the functions $X, Y, P, Q$ on $\B$ would satisfy the \BB/ condition
\begin{equation}\label{bcond}
\{dP \& dX, dQ \& dY\} \equiv \{\pi_1 \& \eta^1, \pi_2 \& \eta^2\}
\mod \theta,\thetabar
\end{equation}
(see the definition at the beginning of \S2).

Accordingly, we let $\B = \scM \times \R$, with coordinate $Z$ on the second factor.
We will show that, on an open neighborhood of any point of $\B$,
 there exist functions $X,Y,P,Q$ such that the ideal
$\J = \I \cup \{\thetabar\}$ on $\B$ (where $\thetabar$ is defined as in \eqref{defthetabar})
gives a \BT/ between $(\scM, \I)$ and $(\Mbar, \Ibar)$.
We will do this by setting up an EDS whose integral manifolds correspond to functions satisfying these conditions;
once we know that this EDS is involutive, an application of the Cartan-K\"ahler Theorem
 will prove the existence of the desired \BT/s.

Let $\B' = \scM' \times \R$, again with $Z$ as the coordinate on the second factor;
we extend the projection $\scM' \to \scM$ to a projection $\B' \to \B$ by the identity
on the second factor.  It will be convenient to set up our EDS in terms of the
coframing $(\theta,\thetabar,\theta_1,\theta_2,\eta^1,\eta^2,\pi_3,\pi_4)$ on $\B'$.
Thus, we will regard $X, Y, P, Q$ as functions on $\B'$,
but require that
\[ dX, dY, dP, dQ \in \{\theta, \thetabar, \theta_1, \theta_2, \eta^1, \eta^2\} \]
so that they are in fact well-defined on $\B$.
In order to satisfy the \BB/ condition \eqref{bcond}, we will furthermore require that
\begin{equation}
 \{dX, dP\} \subset \{\theta, \thetabar, \theta_1, \eta^1\}, \qquad \{dY, dQ\} \subset \{\theta, \thetabar, \theta_2, \eta^2\}.  \label{inclusion1}
\end{equation}
From this, and the structure equations \eqref{Darboux-section-structure},
it follows that $\{dX, dP\}$ (resp., $\{dY, dQ\}$) is the largest integrable subsystem of
$\{\theta, \thetabar, \theta_1, \eta^1\}$ (resp., $\{\theta, \thetabar, \theta_2, \eta^2\}$).  Therefore,
\[ \eta^1 = dx \in \{dX, dP\}, \qquad \eta^2 = dy \in \{dY, dQ\}, \]
and by a contact transformation on $\Mbar$, we may assume that $X=x,\ Y=y$.  Thus, we will set
\[  \thetabar = dZ - P\, dx - Q\, dy, \]
and condition \eqref{inclusion1} becomes
\[ dP \in \{\theta, \thetabar, \theta_1, \eta^1\}, \qquad dQ \in \{\theta, \thetabar, \theta_2, \eta^2\}.  \]
Suppose that
\begin{align}
dP & = P_0 \theta + P_{\overline{0}} \thetabar + P_1 \theta_1 + P_3 \eta^1 \label{nmi-dP-eq-1} \\
dQ & = Q_0 \theta + Q_{\overline{0}} \thetabar + Q_2 \theta_2 + Q_4 \eta^2 .\label{nmi-dQ-eq-1}
\end{align}
Observe that normality of the \BT/ requires that $P_1, Q_2 \neq 0$ and $P_1 \neq Q_2$.

\begin{rem}
Equations \eqref{nmi-dP-eq-1}-\eqref{nmi-dQ-eq-1} give an overdetermined system of first-order
partial differential equations for functions $P$ and $Q$.  
The process of generating compatibility conditions
for such systems can be carried out systematically by computing the exterior derivatives of the 1-form
equations, and using the fact that the repeated exterior derivative of a function is zero.
Moreover, applying Cartan's Test for involutivity (see \cite{CFB}, Chapter 7)
to the resulting EDS will tell us
when we can stop differentiating:  if the system is involutive then
no further compatibility conditions arise through differentiation,
and solutions exist that may be constructed by applying
the Cartan-K\"ahler Theorem.
\end{rem}

Differentiating \eqref{nmi-dP-eq-1} modulo $\theta, \thetabar, \theta_1, \eta^1$ yields
\[ (P_0 + P_{\oo} Q_2)\, \theta_2 \& \eta^2 = 0, \]
and differentiating \eqref{nmi-dQ-eq-1} modulo $\theta, \thetabar, \theta_2, \eta^2$ yields
\[ (Q_0 + Q_{\oo} P_1)\, \theta_1 \& \eta^1 = 0. \]
Therefore, because the 1-forms $\theta, \thetabar, \theta_1, \theta_2, \eta^1, \eta^2$ are linearly independent on $\B$, we have $P_0 = -P_{\oo} Q_2,\ Q_0 = -Q_{\oo} P_1$, and we may write
\begin{align}
dP & = P_{\overline{0}} (\thetabar - Q_2 \theta) + P_1 \theta_1 + P_3 \eta^1 \label{nmi-dPQ-eqs1} \\
dQ & = Q_{\overline{0}} (\thetabar - P_1 \theta) + Q_2 \theta_2 + Q_4 \eta^2 .\notag
\end{align}
Note that neither $P_{\oo}$ nor $Q_{\oo}$ can vanish identically: for, if $P_{\oo} = 0$, then differentiating \eqref{nmi-dPQ-eqs1} shows that $P_1 = 0$ as well, which contradicts the hypothesis of normality.  (A similar argument applies to $Q_{\oo}$.)

Differentiating equations \eqref{nmi-dPQ-eqs1} modulo various combinations of 1-forms leads to the following expressions for the exterior derivatives of $P_{\oo}, P_1, P_3, Q_{\oo}, Q_2, Q_4$:
\begin{align}
dP_{\oo} & = -P_{\oo} Q_{\oo 2} \theta + P_{\oo \oo} (\thetabar - Q_2 \theta) + P_{\oo 1} \theta_1 + P_{\oo 3} \eta^1 \notag \\
dP_1 & = P_{\oo 1} (\thetabar - Q_2 \theta) - \frac{(A_1 Q_2 + C_1 Q_{\oo} P_1)}{Q_{\oo}} \theta_1 - C_2 P_1 \theta_2 \notag \\
& \qquad - \frac{(B_1 Q_2 + Q_{\oo}(E_1 P_1 + P_{\oo}(Q_2 - P_1)))}{Q_{\oo}} \eta^1 - E_2 P_1 \eta^2 + F_1 P_1 \pi_3 + F_2 P_1 \pi_4 \notag \\
dP_3 & = (P_{\oo 3} + P_{\oo}^2) (\thetabar - Q_2 \theta) - \frac{(B_1 Q_2 + Q_{\oo}((D_1 + E_1) P_1 + 2P_{\oo}(Q_2 - P_1)))}{Q_{\oo}} \theta_1 \notag \\
& \qquad + P_{33} \eta^1 + P_1 \pi_3 \label{nmi-dPQ-eqs2} \\
dQ_{\oo} & = -Q_{\oo} P_{\oo 1} \theta + Q_{\oo \oo} (\thetabar - P_1 \theta) + Q_{\oo 2} \theta_2 + Q_{\oo 4} \eta^2 \notag \\
dQ_2 & = Q_{\oo 2} (\thetabar - P_1 \theta) - C_1 Q_2 \theta_1 - \frac{(A_2 P_1 + C_2 P_{\oo} Q_2)}{P_{\oo}} \theta_2 \notag \\
& \qquad - E_1 Q_2 \eta^2 - \frac{(B_2 P_1 + P_{\oo}(E_2 Q_2 + Q_{\oo}(P_1 - Q_2)))}{P_{\oo}} \eta^2  + F_1 Q_2 \pi_3 + F_2 Q_2 \pi_4 \notag \\
dQ_4 & = (Q_{\oo 4} + Q_{\oo}^2) (\thetabar - P_1 \theta) - \frac{(B_2 P_1 + P_{\oo}((D_2 + E_2) Q_2 + 2Q_{\oo}(P_1 - Q_2)))}{P_{\oo}} \theta_2 \notag \\
& \qquad + Q_{44} \eta^2 + Q_2 \pi_4. \notag
\end{align}

Now computing $d(dP_1) \equiv 0$ modulo $\thetabar - Q_2 \theta$ and $d(dQ_2) \equiv 0$ modulo $\thetabar - P_1 \theta$ yields
\begin{align}
P_{\oo \oo} & = \frac{1}{Q_{\oo} P_1^2 Q_2 (Q_2 - P_1)} \Big{(}
2 C_1 (P_1 - 2 Q_2)P_1 Q_2 P_{\oo}^2 Q_{\oo} + 2 C_2 P_1^2 Q_2 P_{\oo} Q_{\oo}^2    \notag \\
& \qquad  + 2 (C_{1,3} - C_1 D_1 + B_1 H) P_1^2 Q_2 P_{\oo} Q_{\oo} +  A_1(P_1 - Q_2) Q_2^2 P_{\oo}^2 + A_2 P_1^3 Q_{\oo}^2 \Big{)} \notag \\
Q_{\oo \oo} & =  \frac{1}{P_{\oo} P_1 Q_2^2 (P_1 - Q_2)} \Big{(}
 2 C_1 P_1 Q_2^2 P_{\oo}^2 Q_{\oo}  + 2 C_2 (Q_2 - 2 P_1)P_1 Q_2 P_{\oo} Q_{\oo}^2 \notag  \\
& \qquad  + 2 (C_{2,4} - C_2 D_2 + B_2 H) P_1 Q_2^2 P_{\oo} Q_{\oo} + A_1 Q_2^3 P_{\oo}^2 +  A_2(Q_2 - P_1) P_1^2 Q_{\oo}^2 \Big{)}\notag  \\
P_{\oo 1} & = -\frac{P_{\oo} (A_1 Q_2 + 2 C_1 Q_{\oo} P_1)}{Q_{\oo} P_1} \label{second-derivatives} \\
Q_{\oo 2} & = -\frac{Q_{\oo} (A_2 P_1 + 2 C_2 P_{\oo} Q_2)}{P_{\oo} Q_2} \notag \\
P_{\oo 3} & = -\frac{P_{\oo}(B_1 Q_2 + (D_1 + E_1) P_1 Q_{\oo}+ (2 Q_2 - P_1) P_{\oo} Q_{\oo})}{P_1 Q_{\oo}} \notag \\
Q_{\oo 4} & = -\frac{Q_{\oo}(B_2 P_1 + (D_2 + E_2) Q_2 P_{\oo}+ (2 P_1 - Q_2) P_{\oo} Q_{\oo})}{Q_2 P_{\oo}}. \notag
\end{align}
This leaves only $P_{33}$ and $Q_{44}$ as undetermined second derivatives of $P$ and $Q$.

We are now ready to set up our exterior differential system.  Let $\widehat{\B} = \B \times \R^{10}$, with coordinates $P$, $Q$, $P_{\oo}$, $P_1$, $P_3$, $Q_{\oo}$, $Q_2$, $Q_4$, $P_{33}$, $Q_{44}$ on the $\R^{10}$ factor.  Let $\W$ be the rank 8 Pfaffian EDS on $\widehat{\B}$ generated by the 1-forms
\begin{align*}
\Theta_1 & = dP - P_{\overline{0}} (\thetabar - Q_2 \theta) - P_1 \theta_1 - P_3 \eta^1, \\
\Theta_2 & = dQ - Q_{\overline{0}} (\thetabar - P_1 \theta) - Q_2 \theta_2 - Q_4 \eta^2,
\end{align*}
and similar forms $\Theta_3, \ldots, \Theta_8$ prescribing conditions \eqref{nmi-dPQ-eqs2} for $dP_{\oo}$, $dP_1$, $dP_3$, $dQ_{\oo}$, $dQ_2$, and $dQ_4$, substituting the values \eqref{second-derivatives} for $P_{\oo \oo}$, $P_{\oo 1}$, $P_{\oo 3}$, $Q_{\oo \oo}$, $Q_{\oo 2}$, $Q_{\oo 4}$.  Integral manifolds of $\W$ satisfying the independence condition
$\theta \& \thetabar \& \theta_1 \& \theta_2 \& \eta^1 \& \eta^2 \& \pi_3 \& \pi_4 \neq 0$
are in one-to-one correspondence with the desired functions $P,Q$ defining a \BT/.

The structure equations for this EDS have the form:
\begin{equation}
\begin{bmatrix} d\Theta_1 \\[0.1in] d\Theta_2 \\[0.1in] d\Theta_3 \\[0.1in] d\Theta_4 \\[0.1in] d\Theta_5 \\[0.1in] d\Theta_6 \\[0.1in] d\Theta_7 \\[0.1in] d\Theta_8 \end{bmatrix} \equiv
- \begin{bmatrix}
0 & 0 & 0 & 0 & 0 & 0 & 0 & 0 \\[0.1in]
0 & 0 & 0 & 0 & 0 & 0 & 0 & 0 \\[0.1in]
0 & 0 & 0 & 0 & 0 & 0 & 0 & 0 \\[0.1in]
0 & 0 & 0 & 0 & 0 & 0 & 0 & 0 \\[0.1in]
0 & 0 & 0 & 0 & \Pi_1 & 0 & 0 & 0 \\[0.1in]
0 & 0 & 0 & 0 & 0 & 0 & 0 & 0 \\[0.1in]
0 & 0 & 0 & 0 & 0 & 0 & 0 & 0 \\[0.1in]
0 & 0 & 0 & 0 & 0 & \Pi_2 & 0 & 0 \\[0.1in]
\end{bmatrix} \&
\begin{bmatrix} \theta \\[0.1in] \thetabar \\[0.1in] \theta_1 \\[0.1in] \theta_2 \\[0.1in] \eta^1 \\[0.1in] \eta^2 \\[0.1in] \pi_3 \\[0.1in] \pi_4 \end{bmatrix} +
\begin{bmatrix} \Psi_1 \\[0.1in] \Psi_2 \\[0.1in] \Psi_3 \\[0.1in] \Psi_4 \\[0.1in] \Psi_5 \\[0.1in] \Psi_6 \\[0.1in] \Psi_7 \\[0.1in] \Psi_8
\end{bmatrix}
\mod{\Theta_1, \ldots, \Theta_8},  \label{EDS-structure-eqs}
\end{equation}
where
\[ \left.
\begin{aligned}
\Pi_1 & \equiv dP_{33} \\
\Pi_2 & \equiv dQ_{44}
\end{aligned} \right\} \mod{\theta, \thetabar, \theta_1, \theta_2, \eta^1, \eta^2, \pi_3, \pi_4},
\]
and $\Psi_1, \ldots, \Psi_8$ are 2-forms which are quadratic in the forms $\theta, \thetabar, \theta_1, \theta_2, \eta^1, \eta^2, \pi_3, \pi_4$, with coefficients which are polynomial functions of the two quantities
\begin{align}
& \big{(} 2H P_1(Q_2 - P_1) Q_{\oo}^2 + Q_2(2 C_1 Q_2+ J_1(P_1 - Q_2)) Q_{\oo} + A_1 Q_2^2 \big{)} P_{\oo} \notag \\
& \qquad \qquad \qquad  - 2 C_2 P_1^2 Q_{\oo}^2 + 2 (C_{2,4} - C_2 D_2 + B_2 H) Q_{\oo} P_1 Q_2 , \label{zero4}
\end{align}
\begin{align}
& \big{(} 2H Q_2(P_1 - Q_2) P_{\oo}^2 + P_1(2 C_2 P_1+ J_2(Q_2 - P_1)) P_{\oo} + A_2 P_1^2 \big{)} Q_{\oo} \notag \\
& \qquad \qquad \qquad  - 2 C_1 Q_2^2 P_{\oo}^2 + 2 (C_{1,3} - C_1 D_1 + B_1 H) P_{\oo} P_1 Q_2 . \label{zero5}
\end{align}
In order to find integral manifolds, we must restrict $\W$ to the locus $\Z \subset \widehat{\B} $ defined by the simultaneous vanishing of \eqref{zero4} and \eqref{zero5}.  (These relations could also be found more directly, by computing $d(dP_1) = d(dQ_2) = 0$.)
While relations of this sort could easily lead to hopeless incompatibility for the PDE system given by \eqref{nmi-dPQ-eqs1} and \eqref{nmi-dPQ-eqs2}, it turns out that differentiating these quantities yields no new relations.

A case-by-case analysis, based on the vanishing or non-vanishing of various torsion coefficients of $(\scM, \I)$, shows that the functions $P_{\oo}, Q_{\oo}, P_1, Q_2$ are all nonzero on an open subset $\Z^0 \subset \Z$ which is surjective for the projection $\Z \to \B$.  Since normality requires that these functions be generically nonzero, we further restrict $\W$ to this open subset.

Let $\overline{\W}$ denote the pullback of $\W$ to $\Z^0$;
$\overline{\W}$ is a rank 6 Pfaffian EDS on $\Z^0$.  Because differentiating the equations defining $\Z$ yields no new relations, $\overline{\W}$ is torsion-free; moreover, it is straightforward to check that $\overline{\W}$ is involutive with last nonzero Cartan character $s_1 = 2$. (See \cite{CFB} for a discussion of Cartan's test and involutivity.)
Therefore, it follows from the Cartan-K\"ahler Theorem that local integral manifolds exist and are parametrized by 2 functions of one variable.

We summarize this result as:

\begin{prop}
Let $(\scM, \I)$ be a hyperbolic Monge-Amp\`ere system which is not Monge-integrable,
and which is Darboux-integrable after one prolongation.  Then near
any point $p \in \scM$ there is an open set $U\subset\scM$ around $p$ such that
the restriction of $\I$ to $U$ is linked to the wave equation by a normal \BT/; moreover, the set of all such \BT/s is parametrized by 2 functions of one variable.
Up to contact transformations, all such B\"acklund transformations
preserve the space of independent variables $x,y$.
\end{prop}

\subsection{The Monge-integrable case} \label{MIcase}\label{DtoBmonge}

In this subsection $(\scM, \I)$ is assumed to be a hyperbolic \MA/ system which is Monge-integrable and Darboux-integrable after one prolongation.  As explained below, we will construct a canonical coframing associated to the partial prolongation of $\I$, and then proceed as in \S \ref{NMIcase}.

A similar argument to that of Lemma \ref{NMIm5coframe} can be used to prove:
\begin{lemma}\label{MIm5coframe} Near any point of $\scM$, there exists a coframing $(\theta,\pi_1,\pi_2,\eta^1,\eta^2)$
such that $\theta$ spans the 1-forms of $\I$, and the characteristic systems
$\calC_1, \calC_2$ of $\I$ have
derived flags
$$
\calC_1 = \{ \theta,\pi_1,\eta^1\} \supset \{\pi_1, \eta^1\} \supset \{\eta^1 \} =\calC_1^{(\infty)},
\qquad
\calC_2 = \{ \theta,\pi_2,\eta^2\} \supset \{\pi_2, \eta^2\} =\calC_2^{(\infty)}.
$$
\end{lemma}
Indeed, the same coframing as that given in the proof of Lemma \ref{NMIm5coframe} satisfies the conditions of Lemma \ref{MIm5coframe}.  Note that this lemma only assumes the Monge-integrability of $(\scM, \I)$.

In terms of the local coframing on $\scM$ given by the lemma, the
partial prolongation $(\scM', \I')$ is defined as follows:
let $\scM' = \scM \times \R$, with coordinate $r$ on the $\R$ factor,
and let $\I'$ be the Pfaffian system on $\scM'$
generated by $\theta$, the 1-form $\theta_1 = \pi_1 - r \eta^1$, and the 2-form $\pi_2 \& \eta^2$.

\begin{lemma}\label{MIn6coframe} Near any point of $\scM'$ there exists a coframing
$(\theta,\theta_1,\pi_2,\eta^1,\eta^2,\pi_3)$ such that
$\I'$ is generated by $\theta,\theta_1$, and $\pi_2 \& \eta^2$, satisfying
\begin{equation}\label{mi-Darboux-structure-0}
\begin{aligned}
d\theta & = -\theta_1 \& \eta^1 - \theta_2 \& \eta^2 \mod{\theta} 
\\
d\theta_1 & = -\pi_3 \& \eta^1 \mod{\theta, \theta_1}, 
\end{aligned}
\end{equation}
with the derived flags of the characteristic systems of $\I'$ given by
\begin{gather*}
\Cprime1  = \{ \theta, \theta_1, \eta^1, \pi_3 \}  \supset
 \{ \theta_1, \eta^1, \pi_3 \} \supset  \{ \eta^1, \pi_3 \}  = \Cprime1^{(\infty)},
\\
\Cprime2  = \{ \theta, \theta_1, \pi_2, \eta^2\} \supset
\{\theta, \pi_2, \eta^2\} \supset  \{\pi_2, \eta^2\} = \Cprime2^{(\infty)}.
\end{gather*}
\end{lemma}

\begin{proof}
As in the proof of Lemma \ref{n7coframe}, the usual or ``full" prolongation
of $\I$ on $\scM \times \R^2$ is generated by $\theta$, $\theta_1$ and
$\theta_2 = \pi_2 - t\, \eta^2$.  We may construct the coframing $(\theta, \theta_1, \theta_2, \eta^1, \eta^2, \pi_3, \pi_4)$ on $\scM \times \R^2$ precisely as in Lemma \ref{n7coframe}, and this coframing satisfies the  structure equations \eqref{nmi-Darboux-structure-0}.

The hypothesis of Darboux-integrability implies that the characteristic system
$$\K_1 = \{ \theta, \theta_1,\theta_2, \eta_1, \pi_3\}$$ of the prolongation contains a rank 2 Frobenius system.  As in Lemma \ref{n7coframe}, it follows from the structure equations and the
construction of $\pi_3$ that
$$\K_1^{(1)} = \{ \theta, \theta_1, \eta_1, \pi_3\}.$$
However, this system is well-defined on $\scM' = \scM \times \R$, and it coincides with
the characteristic system $\Cprime1$ given in the statement of the present Lemma.
It follows that we may adjust $\pi_3$ so that it lies in the rank 2 Frobenius system $\Cprime1^{(\infty)}$.

Note that the second characteristic system $\Cprime2$ is simply the sum of $\I'$ and
the pullback of the characteristic system $\calC_2$ of $\I$, and the structure of its
derived flag follows from that of $\calC_2$.
\end{proof}

The conditions in Lemma \ref{MIn6coframe} are preserved by changes of coframing of the form
\begin{equation} \label{MI-Darboux-freedom}
\begin{bmatrix}
\tilde{\theta}\\[0.1in] \tilde{\theta}_1 \\[0.1in] \tilde{\pi}_2 \\[0.1in] \tilde{\eta}^1 \\[0.1in] \tilde{\eta}^2 \\[0.1in] \tilde{\pi}_3
\end{bmatrix} =
\begin{bmatrix}
c & 0 & 0 & 0 & 0 & 0  \\[0.1in]
0 & a_1 c & 0 & 0 & 0 & 0   \\[0.1in]
0 & 0 & a_2 c & 0 & b_1 & 0   \\[0.1in]
0 & 0 & 0 & a_1^{-1} & 0 & 0   \\[0.1in]
0 & 0 & b_2 & 0 & a_2^{-1} & 0    \\[0.1in]
0 & 0 & 0 & b_3 & 0 & a_1^2 c
\end{bmatrix}^{-1}
\begin{bmatrix}
\theta\\[0.1in] \theta_1 \\[0.1in] \pi_2 \\[0.1in] \eta^1 \\[0.1in] \eta^2 \\[0.1in] \pi_3
\end{bmatrix},
\end{equation}
with $a_1, a_2, c \neq 0$.  Let $G \subset GL(6, \R)$ be the group of such transformations, and let $\calP$ be the $G$-structure on $\scM'$ of which the coframing of Lemma \ref{MIn6coframe} is a section.

After absorbing as much torsion as possible and differentiating to uncover relations among the torsion, $\calP$ has structure equations
\begin{multline}
\begin{bmatrix} d\theta\\[0.1in] d\theta_1\\[0.1in] d\pi_2 \\[0.1in] d\eta^1\\[0.1in] d\eta^2 \\[0.1in] d\pi_3
\end{bmatrix} = -
{\setlength\arraycolsep{0.5pt}
\left[  \begin{array}{c c c c c c}
\gamma & 0 & 0 & 0 & 0 & 0   \\[0.1in]
0 & \ \gamma + \alp_1 & 0 & 0 & 0 & 0  \\[0.1in]
0 & 0 & \gamma + \alp_2 & 0 & \beta_1 & 0  \\[0.1in]
0 & 0 & 0 & -\alp_1 & 0 & 0  \\[0.1in]
0 & 0 & \beta_2 & 0 & -\alp_2 & 0  \\[0.1in]
0 & 0 & 0 & \bet_3 & 0 & \gamma + 2\alp_1
\end{array}\right]
}
 \&
\begin{bmatrix} \theta\\[0.1in] \theta_1\\[0.1in] \pi_2 \\[0.1in] \eta^1\\[0.1in] \eta^2 \\[0.1in] \pi_3
\end{bmatrix}
\label{MI-Darboux-structure}
- \begin{bmatrix}
\theta_1 \& \eta^1 \!  +  \! \pi_2 \& \eta^2 \\[0.1in]
\pi_3 \& \eta^1  \! + \! (A \pi_2 \! + \! B \eta^2) \& \theta \\[0.1in]
0 \\[0.1in]
0 \\[0.1in]
0 \\[0.1in]
C \pi_3 \& \theta_1
\end{bmatrix}.
\end{multline}
Because of the dimensions of the derived flags of the characteristic systems (given in Lemma \ref{MIn6coframe}), $A, B$ are not both zero.  Furthermore, we can choose a local section $\sigma: \scM' \to \calP$ satisfying the conditions that $\eta^1 = dx$, $\eta^2 = dy$, $\pi_2$ is exact, and $\pi_3$ is integrable.
The resulting coframing is uniquely determined.

When we pull back the structure equations via $\sigma$, the pseudoconnection forms $\alp_1$, $\alp_2$, $\bet_1$, $\bet_2$, $\bet_3$, $\gamma$ become semi-basic.  By making use of the remaining ambiguity in these forms and the conditions imposed thus far on the coframing, we can assume that
$$
\alp_1  = E\eta^1, \qquad
\alp_2  = 0,\qquad
\beta_1  = G \pi_2,\qquad
\beta_2  = 0, \qquad
\beta_3  = (-H + 2E) \pi_3, \qquad
\gamma  = F \pi_2 + G \eta^2
$$
for some functions $E,F,G,H$.
Then the structure equations for this coframing become:
\begin{align}
d\theta & = \theta \& (F \pi_2 + G \eta^2)  - \theta_1 \& \eta^1 - \pi_2 \& \eta^2   \notag \\
d\theta_1 & = \theta_1 \& (F \pi_2 + E \eta^1 + G \eta^2)   - \pi_3 \& \eta^1 + \theta \& (A \pi_2 + B \eta^2)   \notag \\
d\pi_2 & = 0 \label{MI-Darboux-section-structure} \\
d\eta^1 & = 0 \notag   \\
d\eta^2 & =  0  \notag \\
d\pi_3 & = \pi_3 \& (-C \theta_1 + F \pi_2 + H \eta^1 + G \eta^2) .    \notag
\end{align}

Once again, we will need to compute relations among the derivatives of the torsion functions in order to show that $(\scM, \I)$ has a \BT/ to the wave equation.  We begin by differentiating the structure equations \eqref{MI-Darboux-section-structure}.
We denote derivatives as, e.g.,
\[ dA = A_{0} \theta + A_{1} \theta_1 + A_{2} \pi_2 + A_{3} \eta^1 + A_{4} \eta^2 + A_{5} \pi_3. \]

Computing $d(d\theta) = d(d\theta_1) = d(d\pi_3) = 0$ yields the following equations for the derivatives of the torsion functions:
\begin{align}
dA & = A_{0} \theta + AC \theta_1 + A_2 \pi_2 - AE \eta^1 + A_4 \eta^2  \notag \\
dB & = B_{0} \theta + BC \theta_1 + A_4 \pi_2 - BE \eta^1 + B_4 \eta^2  \notag \\
dC & = C_1 \theta_1 + CF \pi_2 + C_3 \eta^1 + CG \eta^2 + C_5 \pi_3
\notag \\
dE & = E_1 \theta_1 + 2A \pi_2 + E_3 \eta^1 + 2B \eta^2 - C \pi_3
\label{MI-dtorsion} \\
dF & = AC \theta + F_2 \pi_2 + A \eta^1 + F_4 \eta^2
\notag \\
dG & = BC \theta + F_4 \pi_2 + B \eta^1 + G_4 \eta^2
\notag \\
dH & = (CE - C_3) \theta_1 + A \pi_2 + H_3 \eta^1 + B \eta^2 + H_5 \pi_3 .
\notag
\end{align}

We may obtain further relations among the derivatives of the torsion functions by differentiating equations \eqref{MI-dtorsion}.  Computing $d(dF) \equiv d(dG) \equiv 0$ modulo $\pi_2, \eta^2$ and recalling that $A,B$ cannot vanish simultaneously yields
$$
A_0 = A(C_3 - CE), \qquad
B_0 = B(C_3 - CE), \qquad
C_1 = -C^2, \qquad
C_5 = 0.
$$
Then computing $d(dA) \equiv d(dB) \equiv  0$ modulo $\theta, \pi_2, \eta^2$ yields
\[ E_1 = 2(CE - C_3), \]
and $d(dC)=0$ implies that
\[ dC_3 = (C^2E - 2C C_3) \theta_1 + (AC + F C_3) \pi_2 + C_{33} \eta^1 + (BC + G C_3) \eta^2 - C^2 \pi_3. \]
Finally, computing
$d(dA) \equiv d(dB) \equiv  0$ modulo $\pi_2, \eta^2$ yields
\[ C_{33} = C E_3 + E C_3. \]
We now have all the relations that will be needed for the involutivity calculation below.

Now suppose that $\Mbar=\R^5$ carries a \MA/ system $\Ibar$ representing
the wave equation $Z_{XY}=0$, generated algebraically by the contact form
\begin{equation}\label{MI-defthetabar}
\thetabar = dZ - P\, dX - Q\, dY
\end{equation}
and the 2-forms $dP \& dX$ and $dQ \& dY$.  As in \S \ref{NMIcase},
if there were a B\"acklund transformation $\B \subset \scM \times \Mbar$,
then $Z$ would be a local coordinate on the fibers of $\B \to \scM$
and the functions $X, Y, P, Q$ on $\B$ would satisfy the B\"acklund condition
\begin{equation}\label{MI-bcond}
\{dP \& dX, dQ \& dY\} \equiv \{\pi_1 \& \eta^1, \pi_2 \& \eta^2\}
\mod \theta,\thetabar .
\end{equation}

As in \S \ref{NMIcase}, let $\B' = \scM' \times \R$, again with $Z$ as the coordinate on the second factor;
we extend the projection $\scM' \to \scM$ to a projection $\B' \to \B$ by the identity
on the second factor.  We will regard $X, Y, P, Q$ as functions on $\B'$,
but require that
\[ dX, dY, dP, dQ \in \{\theta, \thetabar, \theta_1, \pi_2, \eta^1, \eta^2\} \]
so that they are in fact well-defined on $\B$.
In order to satisfy the \BB/ condition \eqref{MI-bcond}, we will furthermore require that
\begin{equation}
 \{dX, dP\} \subset \{\theta, \thetabar, \theta_1, \eta^1\}, \qquad \{dY, dQ\} \subset \{\theta, \thetabar, \pi_2, \eta^2\}.  \label{MI-inclusion1}
\end{equation}

The same argument as that given in \S \ref{NMIcase} shows that by a contact transformation on $\Mbar$, we may assume that $X=x$.  However, the same is not true for $Y$: the system $\{\theta, \thetabar, \pi_2, \eta^2\}$ on $\B$ contains a rank 3 integrable subsystem, so we cannot necessarily arrange to have $\eta^2 \in \{dY, dQ\}$.  There are three different, geometrically natural conditions that we could impose on the intersection of the rank 2 Pfaffian systems $\{dY, dQ\}$ and $\{\pi_2, \eta^2\}$, each of them potentially leading to a different type of \BT/:
\begin{enumerate}
\item $\{dY, dQ\} \cap \{\pi_2, \eta^2\}$ has rank 1 and is spanned by a non-integrable 1-form.
\item $\{dY, dQ\} \cap \{\pi_2, \eta^2\}$ has rank 1 and is spanned by an integrable 1-form.
\item $\{dY, dQ\} \cap \{\pi_2, \eta^2\}$ has rank 2.
\end{enumerate}
In cases (2) and (3) we can arrange that $Y=y$ via contact transformations on $\scM$ and $\Mbar$, but in case (1) this is not possible.

\subsubsection{Case (1)}
In this case we have
\[ \thetabar = dZ - P\, dx - Q\, dY, \]
and condition \eqref{MI-inclusion1} becomes
\[ dP \in \{\theta, \thetabar, \theta_1, \eta^1\}, \qquad dY, dQ \in \{\theta, \thetabar, \pi_2, \eta^2\}.  \]
Suppose that
\begin{align}
dP &= P_0 \theta + P_{\overline{0}} \thetabar + P_1 \theta_1 + P_3 \eta^1 \notag \\
dQ &= Q_0 \theta + Q_{\overline{0}} \thetabar + Q_2 \pi_2 + Q_4 \eta^2 \label{mi-dP-eqs-case1} \\
dY &= Y_0 \theta + Y_{\overline{0}} \thetabar + Y_2 \pi_2 + Y_4 \eta^2. \notag
\end{align}
Normality of the \BT/ requires that $P_1 \neq 0, \ Q_2 Y_4 - Q_4 Y_2 \neq 0$, and $P_1 \neq Q_2 Y_4 - Q_4 Y_2.$

The argument proceeds in much the same fashion as that of \S \ref{NMIcase}: differentiating equations \eqref{mi-dP-eqs-case1} leads to relations among the derivatives of $P,Q,Y$.  Eventually we are led to a Pfaffian exterior differential system $\W$ whose integral manifolds satisfying the independence condition $\theta \& \thetabar \& \theta_1 \& \pi_2 \& \eta^1 \& \eta^2 \neq 0$ are in one-to-one correspondence with the desired functions $P,Q,Y$ defining a \BT/.  This EDS is involutive with last nonzero Cartan character $s_3 = 1$.
Therefore, local integral manifolds exist and are parametrized by 1 function of three variables.

If we impose the additional condition that the \BT/ be holonomic, we find that the resulting EDS is involutive with last nonzero Cartan character $s_2 = 2$.  Therefore, among the \BT/s of this type, there is a small, proper subset, parametrized by 2 functions of two variables, consisting of holonomic transformations.

\subsubsection{Case (2)}
In this case, we can use contact transformations on $\scM$ and $\Mbar$ to arrange that $\{dY, dQ\} \cap \{\pi_2, \eta^2\}$ is spanned by $\eta^2 = dy = dY$.  Then we have
\[ \thetabar = dZ - P\, dx - Q\, dy, \]
and condition \eqref{MI-inclusion1} becomes
\[ dP \in \{\theta, \thetabar, \theta_1, \eta^1\}, \qquad dQ \in \{\theta, \thetabar, \pi_2, \eta^2\}.  \]
Suppose that
\begin{align}
dP &= P_0 \theta + P_{\overline{0}} \thetabar + P_1 \theta_1 + P_3 \eta^1, \label{mi-dP-eqs-case2} \\
dQ &= Q_0 \theta + Q_{\overline{0}} \thetabar + Q_2 \pi_2 + Q_4 \eta^2. \notag
\end{align}
Normality of the \BT/ requires that $P_1, Q_2 \neq 0$ and $P_1 \neq Q_2.$

Differentiating equations \eqref{mi-dP-eqs-case2} leads to relations among the derivatives of $P,Q$, and to a Pfaffian exterior differential system $\W$ whose integral manifolds satisfying the independence condition $\theta \& \thetabar \& \theta_1 \& \pi_2 \& \eta^1 \& \eta^2 \neq 0$ are in one-to-one correspondence with the desired functions $P,Q$ defining a \BT/.  This EDS is involutive with last nonzero Cartan character $s_2 = 1$.
Therefore, local integral manifolds exist and are parametrized by 1 function of two variables.

If we impose the additional condition that the \BT/ be holonomic, we find that the resulting EDS is involutive with last nonzero Cartan character $s_1 = 3$.  Therefore, among the \BT/s of this type, there is a small, proper subset, parametrized by 3 functions of one variable, consisting of holonomic transformations.

\subsubsection{Case (3)}
In this case, we can use contact transformations on $\scM$ and $\Mbar$ to arrange that $\eta^2 = dy = dY, \ \pi_2 = dQ$.  Then we have
\[ \thetabar = dZ - P\, dx - Q\, dy, \]
and condition \eqref{MI-inclusion1} becomes
\[ dP \in \{\theta, \thetabar, \theta_1, \eta^1\}.  \]
Suppose that
\begin{equation}
dP = P_0 \theta + P_{\overline{0}} \thetabar + P_1 \theta_1 + P_3 \eta^1. \label{mi-dP-eqs-case3}
\end{equation}
Normality of the \BT/ requires that $P_1 \neq 0$ and $P_1 \neq 1$.

Differentiating equation \eqref{mi-dP-eqs-case3} leads to relations among the derivatives of $P$, and to a Pfaffian exterior differential system $\W$ whose integral manifolds satisfying the independence condition $\theta \& \thetabar \& \theta_1 \& \pi_2 \& \eta^1 \& \eta^2 \neq 0$ are in one-to-one correspondence with the desired functions $P$ defining a \BT/.

The involutivity calculation in this case depends on the torsion functions in the structure equations \eqref{MI-Darboux-section-structure}.  If
\begin{equation}
 AG - BF = AC + F(CE - C_3) = AF^2 + F A_2 - A F_2 = BF^2 + F A_4 - A F_4 = 0, \label{mi-case3-cond}
\end{equation}
then $\W$ is involutive with last nonzero Cartan character $s_2 = 1$, and so local integral manifolds exist and are parametrized by 1 function of two variables.   Otherwise, there are no solutions with $P_1 \neq 0$, and hence no normal \BT/s of this type.

Observe that in this case, the $G$-structure on the B\"acklund transformation $\B$ (cf. \S \ref{wsetup}) will satisfy the condition that (omitting obvious pullback notations)
\[ \{\om^3, \om^4\} = \{\pi_2, \eta^2\} = \{dY, dQ\}. \]
Therefore, all transformations of this type satisfy the hypotheses of Proposition \ref{one-C-vector-zero-prop} and so are holonomic.

We summarize these results as:

\begin{prop}
Let $(\scM, \I)$ be a hyperbolic Monge-Amp\`ere system which is Monge-integrable, and Darboux-integrable after one prolongation.  Then there exist B\"acklund transformations  of types (1) and (2) above between $(\scM, \I)$ and the standard wave equation $(\Mbar, \Ibar)$, and of type (3) if the torsion functions of $(\scM, \I)$ satisfy \eqref{mi-case3-cond}.
The generic B\"acklund transformation is of type (1) and does not preserve the space of independent variables.  There are both holonomic and non-holonomic \BT/s of types (1) and (2), and all B\"acklund transformations of type (3) are holonomic.
\end{prop}

\section{Examples}\label{examples}
\def\Fsigma{F}
\def\wavephi{\thetabar}
In this section we review the classifications of second-order Darboux-integrable \MA/ equations,
due to Goursat and Vessiot, and discuss the connection between our results and
the work of Zvyagin.  We will also give examples of a method for explicitly solving for \BT/s linking
these equations to the wave equation.

\subsection{The Goursat-Vessiot List}
Goursat \cite{G1899} studied non-linear PDE of the form
\begin{equation}\label{sform}
u_{xy} = \Fsigma(x,y,u,u_x,u_y)
\end{equation}
which are Darboux-integrable at the 2-jet level, classifying them up to complex contact transformations that preserve the
form \eqref{sform}.
Using Lie-theoretic techniques,
Vessiot \cite{V1942} reproduced Goursat's classification,
expanded to include linear equations, and showed that some of the equations on Goursat's list were equivalent under more
general contact transformations.  Recently, Biesecker \cite{Mattsthesis} re-proved Vessiot's classification using Cartan's method of equivalence, with respect to real contact transformations.
Retaining Goursat's numbering, the list is:
\begin{align}
(x+y) u_{xy} &= 2\sqrt{u_x u_y}; \tag{I}\\
u\, u_{xy} &= \sqrt{1+u_x^2}\sqrt{1+u_y^2}; \tag{II}\\
(\sin u) u_{xy} &= \sqrt{1+u_x^2}\sqrt{1+u_y^2}; \tag{III} \\
u\, u_{xy} &= \pm \phi(u_x) \psi(u_y), \tag{IV} \\
\intertext{where $\phi(t),\psi(t)$ satisfy the ODE $df/dt \pm t/f= K$ for
some nonzero constant $K$;}
(x+y) u_{xy} &= \gamma(u_x) \gamma(u_y), \tag{V} \\
\intertext{where $\gamma$ is implicitly defined by $\gamma(t)-1 =\exp(t - \gamma(t))$;}
u_x - u\dfrac{u_{xy}}{u_y} &= f\left(x,\dfrac{u_{xy}}{u_y}\right);\tag{VI}\\
u_{xy} &= e^u \sqrt{1+(u_x)^2};\tag{VII} \\
u_x - y\, u_{xy} &= f(x,u_{xy}); \tag{VIII}\\
u_{xy} &= e^u; \tag{IX} \\
u_{xy} &= u_x e^u,\tag{X};\\
u_{xy} &= \left(\dfrac1{u+x}+\dfrac1{u+y}\right)u_x u_y  \tag{XI}. \\
\intertext{(In (VI) and (VIII) the function $f$ is arbitrary.)
To Goursat's original list, Vessiot added representatives of the two equivalence classes of
Darboux-integrable linear equations:}
u_{xy} &= a(x,y)u_x + b(x,y)u_y -a(x,y)b(x,y)u, \tag{XII} \\
\intertext{where $h(x,y)=-a_x$ and $k(x,y)=-b_y$ must satisfy the system $(\ln h)_{xy} = 2h-k$, $(\ln k)_{xy} = 2k-h$
with $h\ne k$; and finally,}
u_{xy} &= \dfrac{2u}{(x+y)^2}.\tag{XIII}
\end{align}

In the above list, we have replaced Goursat's original versions of (VII) and (XI) by
simpler equations that Vessiot showed were equivalent to them by contact transformations;
see \cite{V1942}, part 2, pages 5 and 6, respectively.
Vessiot also observed that
(VI) is contact-equivalent to (X), and (VIII) may be reduced by a contact transformation
to the special case
\begin{equation}
u_{xy} = \dfrac{u_x}{x+y}. \qquad \tag{VIII*}
\end{equation}

\bigskip
Our Theorem \ref{ourtheorem}, together with Goursat's classification, implies the following
\begin{cor}\label{Cor-3} If a second-order \MA/ PDE for one function of two variables is
linked to the standard wave equation by a normal \BT/ with 1-dimensional fibers,
then the PDE is either equivalent to the wave equation by a contact transformation,
or equivalent to one of the equations (I)-(XIII) in the above list.
\end{cor}

\subsection{Zvyagin's List}\label{zsec}
Zvyagin \cite{Zvyagin2} investigated second-order \MA/ equations linked to the standard
wave equation by a \BT/, and asserted that all such transformations
that are non-holonomic are exhausted by a list of six examples in addition to Liouville's equation.
Zvyagin did not publish a proof of this classification, and did not give
explicit forms for the \MA/ equations for some of the transformations on
his list.   He did give an explicit transformation for Goursat-Vessiot equation (I):
\begin{equation}
\sqrt{p}-\sqrt{P} = \sqrt{\dfrac{Z-u}{x-y}},\quad
\sqrt{q}-\sqrt{Q}=-\sqrt{\dfrac{Z-u}{x-y}},\tag{Z.I}
\end{equation}
where, as in \S4, $Z$ is the solution to the wave equation, with
$x$- and $y$- derivatives $P$ and $Q$.  (The $x$- and $y$-coordinates
are preserved by the transformation.)

Corollary \ref{Cor-3} implies that every one of Zvyagin's
transformations must be identifiable with an equation on the Goursat-Vessiot
list.   We have calculated explicit forms for certain transformations
on Zvyagin's list, and we can identify the following
transformations as belonging to equations (II), (III), and (VII), respectively:
\begin{gather}
\label{G2bt}
Z P - u p = \sqrt{Z^2 - u^2} \sqrt{1+P^2},\quad
Z Q - u q = \sqrt{Z^2 - u^2}\sqrt{1+Q^2};
\tag{Z.II} \\[0.15in]
p = \dfrac{(\sinh(z) +\tfrac12 e^{-w})P + e^{(Z-w)/2}\sqrt{P^2-1}}{-\sin u}, \quad
q = \dfrac{(\sinh(z) -\tfrac12 e^{-w})Q + e^{-(w+Z)/2}\sqrt{Q^2-1}}{-\sin u},
\tag{Z.III}\\
\intertext{where $w$ is related to $u$ and $Z$ by
$\cos u = \cosh Z - \tfrac12 e^{-w}$;}
p = (1-2e^{u+Z})P - 2 e^{(u+Z)/2}\sqrt{e^{u+Z}-1}\sqrt{P^2+1},
\qquad
q = -Q - e^{(u-Z)/2}\sqrt{e^{u+Z}-1}.
\tag{Z.VII}
\end{gather}
Each of the above transformations, which preserve the $x$- and $y$-coordinates,
 may be verified as being non-holonomic.
To see how this is done, suppose that the transformation equations, when solved for $p$ and $q$, take the form
\begin{equation}\label{transeq}
p = f(x,y,u,Z,P), \qquad q = g(x,y,u,Z,Q).
\end{equation}
The Cartan system
for the \MA/ equation is spanned
by $dx,dy,du,dp$ and $dq$, while the Cartan system for the
wave equation is spanned by $dx,dy,dZ, dP,dQ$.
Recall from \S2 that a transformation is holonomic if the intersection
of these systems is Frobenius.  In light of the
transformation equations \eqref{transeq}, the intersection of
these two systems is spanned by $dx$, $dy$, $dp - f_u\,du$
and $dq - g_u\,du$.  The last two 1-forms are congruent modulo $dx$ and $dy$
to
$$\xi_1 =  f_P dP + f_Z dZ, \qquad \xi_2 = g_Q dQ + g_Z dZ,$$
respectively.  To check that the system $\{ dx, dy, \xi_1, \xi_2\}$ is not Frobenius,
compute
\begin{equation}
\label{holonomicxi}
\begin{aligned}
d\xi_1 \wedge \xi_1 \wedge \xi_2 &\equiv g_Q (f_{uZ} f_P - f_{Pu} f_Z)dP \& dQ \& du \& dZ, \\
d\xi_2 \wedge \xi_1 \wedge \xi_2 &\equiv f_P (g_{uZ} g_Q - g_{Qu} g_Z)dP \& dQ \& du \& dZ
\end{aligned}
\mod dx, dy.
\end{equation}
In each case, the coefficients on the right are nonzero, and we
conclude that the transformation is non-holonomic.

\subsection{Solving for \BT/s}\label{solving}
In this subsection, we will set up systems of PDE
whose solutions are \BT/s to the wave equation
for some examples on the Goursat-Vessiot list.
Although the existence of these transformations follows
from the arguments of \S4, here we will be
able to go further in writing down explicit formulas
for the transformations.   Because we will work with
specific \MA/ equations on the list, we can take advantage
of explicit formulas for the characteristic invariants.
(These invariants are computed, for example, in the dissertation
of M. Biesecker \cite{Mattsthesis}.)

The general approach is as follows.  We write a PDE on the list in the form
\begin{equation}\label{sform4}
s=\Fsigma(x,y,u,p,q).
\end{equation}
This form always has $x$ and $y$ as characteristic invariants, and we assume
these are the only functionally independent invariants up to first order
for the equation (i.e., we assume that the equation is not Monge-integrable).
The \BT/ must take these invariants to
corresponding characteristic invariants for the wave equation. By employing
a change of variables on the wave equation side, of the form $X \mapsto \phi(X)$, $Y \mapsto \psi(Y)$, and
interchanging $X$ and $Y$ if necessary, we may assume that the transformation has
$$x=X, \qquad y=Y.$$
Now suppose that the remaining equations defining the \BT/ take the form
\begin{equation}\label{pandqgeneral}
p = f(x,y,u,Z,P,Q), \qquad q = g(x,y,u,Z,P,Q).
\end{equation}

The \MA/ system on $\R^5$ encoding the PDE \eqref{sform4} is generated
algebraically by the contact form $\theta = du - p\, dx - q\, dy$ and the 2-forms
$$\Omega_1 = (dp - \Fsigma(x,y,u,p,q) dy) \& dx,
\qquad \Omega_2 = (dq - \Fsigma(x,y,u,p,q) dx) \& dy.$$
The defining property of the \BT/ is that substituting \eqref{pandqgeneral}
into $\Omega_1, \Omega_2$ must make them congruent to linear combinations
of $dP \& dx$ and $dQ \& dy$ (the 2-forms defining the \MA/ system
for the wave equation) modulo $\theta$ and the contact form on the
wave equation side,
$$\wavephi = dZ - P\,dx - Q\,dy.$$
In fact, $\Omega_1$ must become congruent to a multiple of $dP \& dx$ and
$\Omega_2$ congruent to a multiple of $dQ \& dy$.  Using \eqref{pandqgeneral},
we compute
\begin{align*}
\Omega_1 &\equiv \left((f_y + f_u g + f_Z Q-\Fsigma) dy + f_P dP + f_Q dQ\right) \& dx
\qquad \mod \theta, \wavephi,\\
\intertext{and}
\Omega_2 &\equiv \left((g_x + g_u f + g_Z P - \Fsigma) dx + g_P dP + g_Q dQ\right) \& dy
\qquad \mod \theta, \wavephi.
\end{align*}
We immediately conclude that $f_Q = g_P = 0$, so that
the transformation is of the form
\begin{equation}\label{pandq}
p = f(x,y,u,Z,P), \qquad q = g(x,y,u,Z,Q),
\end{equation}
and $f,g$ must satisfy
two additional first-order PDEs,
\begin{align}
f_y &=\Fsigma(x,y,u,f,g)-  f_u g - f_Z Q, \label{genfy}\\
g_x &=\Fsigma(x,y,u,f,g) - g_u f - g_Z P. \label{gengx}
\end{align}

We derive additional first- and second-order PDEs that $f$ and $g$ must satisfy
by differentiating the conditions so far.
Taking derivatives with respect to $Q$ in \eqref{genfy} and $P$ in \eqref{gengx}
gives
\begin{equation}\label{genfzgz}
f_Z = \left(\Fsigma_q - f_u\right) g_Q,\qquad
g_Z = \left(\Fsigma_p - g_u\right) f_P,
\end{equation}
where the partials $\Fsigma_p = \di\Fsigma/\di p$ and $\Fsigma_q=\di\Fsigma/\di q$ are taken
and then evaluated with $p$ and $q$ given by \eqref{pandq}.   As we will
see in specific cases below, this will sometimes imply that $f$ and $g$
must be linear in $P$ and $Q$.

In what follows, let $J_1$ and $J_2$ denote the second-order characteristic
invariants for the given PDE (whose existence makes the equation
Darboux-integrable), expressed in
terms of $x,y,u,p,q$ and the second-order jet coordinates
$r$ and $t$.  (We make the convention that
$J_1$ is invariant along the characteristic curves where $x$ is constant,
and $J_2$ is invariant when $y$ is constant.)
Then the \BT/ must take $J_1$ and $J_2$
to second-order characteristic invariants for the wave equation.
In order to compute these additional constraints, we must take total
$x$- and $y$-derivatives in \eqref{pandq}
to deduce how the second-order jet coordinates $r$ and $t$
transform in terms of those of the wave equation:
\begin{equation}\label{randt}
r = f_x + f_u f + f_Z P + f_P R, \qquad
t = g_y + g_u g + g_Z Q + g_Q T.
\end{equation}
Requiring that, under these substitutions, $J_1$ transforms to be a function of only  $x,P,R$,
and $J_2$ transforms to be a function of only $y,Q,T$,
will lead to
additional second-order PDEs which $f$ and $g$ must satisfy.

We now turn to specific examples.
%
%

\subsubsection*{Equation IX (Liouville's equation)}
In this case, $\Fsigma = e^u$, and the equations
\eqref{genfy} through \eqref{genfzgz} become
\begin{align}
f_y &= e^u +f_u(Q g_Q - g), \label{Lioufy}\\
g_x &= e^u + g_u( P f_P - f), \label{Liougx}\\
f_Z &= -f_u g_Q,\label{Lioufz} \\
g_Z &= -g_u f_P. \label{Liougz}
\end{align}
As mentioned in \S1, the characteristic invariants are
$$J_1 = r -\tfrac12 p^2, \qquad J_2 =t - \tfrac12 q^2.$$
Under \eqref{pandq}, the first invariant
transforms as
$$r-\tfrac12 p^2 = f_x + f_u f + f_Z P + f_P R - \tfrac12 f^2.$$
Requiring that this be a function of $x,P,R$ only immediately implies
that $f_P$ can depend on $x$ and $P$ only, and that the remaining
terms have no dependence on $u,Z$ or $y$.  This gives us 6
additional second-order PDEs for $f$:
\begin{align}
\label{Liousecond1}
f_{Pu}=f_{Py}=f_{PZ}&=0, \quad & \di_u, \di_y, \di_Z(f_x +f_u f + f_Z P  - \tfrac12 f^2) &=0. \\
\intertext{Similarly, we also get}
\label{Liousecond2}
g_{Qu}=g_{Qy}=g_{QZ}&=0, \quad & \di_u, \di_x, \di_Z (g_y +g_u g + g_Z Q  - \tfrac12 g^2) &=0.
\end{align}
Note that some of these second-order equations are redundant,
in light of the derivatives of \eqref{Lioufz} and \eqref{Liougz}.

Next, we derive additional equations by differentiation.
Note that \eqref{Lioufy} shows that $f_u$ cannot be identically zero;
then, taking a $Q$ derivative of \eqref{Lioufz} shows that  $g_{QQ}=0$.
Similarly, $f_{PP}=0$, so that $f$ and $g$ are linear in $P$ and $Q$.  Thus, we may set
$$f(x,y,u,P,Z) = f^0(x)P + f^1(x,y,u,Z), \qquad
g(x,y,u,Q,Z) = g^0(y)Q + g^1(x,y,u,Z).
$$
In particular, taking the terms in \eqref{Liousecond1}, \eqref{Liousecond2} that are linear
in $P$ and $Q$ respectively gives
\begin{equation}\label{Lioufugu}
\di_u\left((f^0-g^0) f^1_u - f^0 f^1\right) = 0,
\qquad
\di_u\left((g^0-f^0) g^1_u - g^0 g^1\right) = 0.
\end{equation}
Furthermore, equating the $Z$-derivative of \eqref{Lioufy}
with the $y$-derivative of \eqref{Lioufz}, and using the $u$-derivatives of these equations to
determine $f_{yu}$ and $f_{Zu}$, gives the compatibility condition
\begin{align*}
(g^0-f^0) f^1_u g^1_u &= g^0_y f^1_u + e^u g^0;\\
\intertext{we similarly derive}
(f^0-g^0) f^1_u g^1_u &= f^0_x g^1_u + e^u f^0.
\end{align*}
Adding and differentiating with respect to $u$, and using the values
for $f^1_{uu}$ and $g^1_{uu}$ given by \eqref{Lioufugu},
shows that $f^0=-g^0=k$ for some nonzero constant $k$, and $f^1_u g^1_u = \tfrac12 e^u$.
Integrating the remaining equations shows that the most general form for the
transformation is
\begin{equation}\label{Liouvillegeneral}
\begin{aligned}
p &= k P + 2\exp\left(\dfrac{u+k Z + v(x) + w(y)}{2}\right)+v'(x),\\
q &= -k Q + \exp\left(\dfrac{u-k Z - v(x)-w(y)}{2}\right) - w'(y),
\end{aligned}
\end{equation}
where $v(x),w(y)$ are arbitrary functions.

Using the calculation \eqref{holonomicxi}, it is easy to verify that none
of these transformations is holonomic.

\bigskip
In the next two examples, we will analyze the system of PDEs that $f$ and $g$ must satisfy
using the techniques of exterior differential systems.

\bigskip
\subsubsection*{Equation XIII}
This PDE,
$$s = \Fsigma(u,x,y) := \dfrac{2u}{(x+y)^2},$$
has second-order characteristic invariants
$$J_1 =r + \dfrac{2p}{x+y}, \qquad J_2 = t + \dfrac{2q}{x+y}$$
in the $x$- and $y$-directions respectively
(see \cite{Mattsthesis}, Appendix A).
Substituting for $p$ and $r$ from \eqref{pandq} and \eqref{randt}
yields
$$J_1 = \left(f_x+ f_u f + f_Z P+\dfrac{2f}{x+y}\right) + f_P R,$$
so that $f_P$ and the expression in parentheses must be functions of $x$ and $P$ only.
Similarly, we have
$$J_2 = \left(g_y + g_u g + g_Z Q+\dfrac{2g}{x+y}\right) + g_Q T,$$
hence $g_Q$ and the expression in parentheses must be functions of $y$ and $Q$ only.

In this case, \eqref{genfy} through \eqref{genfzgz} specialize to
\begin{align}
f_y &= \Fsigma - (g - Q g_Q)f_u, \label{b5fy} \\
g_x &= \Fsigma -(f - P f_P) g_u, \label{b5gx} \\
f_Z &= -g_Q f_u, \label{b5fz} \\
g_Z &= -f_P g_u. \label{b5gz}
\end{align}
If $f_u$ were identically zero, then $f_Z$ would also be identically zero, but
then $f_y = \Fsigma= 2u/(x+y)^2$ would give a contradiction.  So, we may assume that
$f_u$ and (similarly) $g_u$ are nonzero on an open dense set.
It then follows from \eqref{b5fz} that $g_{QQ}=0$ and from \eqref{b5gz} that $f_{PP}=0$,
i.e., $f$ and $g$ are again linear in $P$ and $Q$.

Differentiating \eqref{b5gx},\eqref{b5gz} with respect to $x$ and $Z$, and equating mixed partials, enables us
to solve for $f_{Px}$ as
\begin{equation}\label{b5fpx}
f_{Px} = f_u (f_P - g_Q) - \dfrac{2 f_P}{g_u(x+y)^2},
\end{equation}
while from \eqref{b5fy},\eqref{b5fz} we similarly obtain
\begin{equation}\label{b5gqy}
g_{Qy} = g_u(g_Q - f_P) - \dfrac{2 g_Q}{f_u(x+y)^2}.
\end{equation}

To encode the PDEs that $f$ and $g$ must satisfy as an exterior differential system,
we will use $x,y,u,Z,P,Q$ as independent variables, and use $f,f_x,f_u,f_P,r_1$
and $g,g_y,g_u,g_Q,t_1$ as dependent variables.
(The role of the coefficients $r_1$ and $t_1$ will
be made clear below.)  We will regard these variables as coordinates on $\R^{16}$.
As stated above, we restrict to the
open subset $\scrU\subset \R^{16}$ where $f_P,f_u,g_Q$ and $g_u$ are nonzero.

The generator 1-forms are $\psi_1$ through $\psi_6$,
where
\begin{align*}
\psi_1 &= -df + f_x dx + f_u du + f_P dP
+ f_Z dZ + f_y dy,
\\
\psi_2 &= -dg + g_y dy + g_u du + g_Q dQ
+ g_Z dZ + g_x dx,
\\
\psi_3 &= -df_P + f_{Px} dx, \\
\psi_4 &= -dg_Q + g_{Qy} dy,
\end{align*}
with $f_y,f_Z,g_x,g_Z$ given by equations \eqref{b5fz} through \eqref{b5gx} and
$f_{Px},g_{Qy}$ given by \eqref{b5fpx} and \eqref{b5gqy}.
The remaining generators $\psi_5,\psi_6$
encode the rest of the condition that the second-order characteristic
invariants be preserved.  Differentiating the first term in $J_1$ gives
\begin{multline*}
d\left(f_x + f_z P + f_u f + \dfrac{2f}{x+y}\right) \equiv
d(f_x) + P d(f_Z) + \left(d(f_u)-\dfrac{2\,dy}{(x+y)^2}\right) f\\
+ \left(f_u+\dfrac{2}{x+y}\right)
\left(f_u \theta -g_Q f_u \wavephi + \Fsigma dy\right)
\mod \psi_1, dx, dP.
\end{multline*}
Let $\eta_1$ be the 1-form on the right; then for any \BT/ $\eta_1$ must
be a linear combination of $dx$ and $dP$.  In fact, since only the first
term in $\eta_1$ can contain $dP$, the coefficient of $dP$ in $\eta_1$ must be $f_{Px}$.  Thus,
our remaining generators are
$$
\psi_5 = \eta_1 - f_{Px} dP - r_1 dx,\qquad
\psi_6 = \eta_2 - g_{Qy} dQ - t_1 dy,
$$
where, based a similar calculation of $dJ_2$, we set
$$\eta_2 = d(g_y) + Q d(g_Z) + \left(d(g_u)-\dfrac{2\,dx}{(x+y)^2}\right)g
+ \left(g_u+\dfrac{2}{x+y}\right)
\left(g_u \theta -f_P g_u \wavephi + \Fsigma dx\right).
$$

We seek to construct integral manifolds of the given
differential ideal, i.e., submanifolds of $\scrU$ to which
all the forms in the ideal pull back to be zero.
An {\em integral element} for an EDS is an infinitesimal
version of an integral manifold, i.e., a subspace in
the tangent space to $\scrU$ at some point, to which all
the forms in the ideal restrict to be zero.  Because
we want integral manifolds which are graphs of functions
of $x,y,u,P,Q,Z$, we will only consider integral elements
which are 6-dimensional, and to which the differentials
$dx,dy,du,dP,dQ,dZ$ restrict to be linearly independent;
we will call these {\em admissible} integral elements.

Applying Cartan's Test to the Pfaffian system
generated by $\psi_1, \ldots, \psi_6$ shows that
it has last nonzero Cartan character $s_1 =4$, but is not involutive,
as the space of admissible integral elements
has 2-dimensional fiber at each point.  However, the system becomes
involutive after one prolongation, and this establishes the existence of
the required \BT/s.  The last nonzero Cartan character of the involutive
prolongation is $s_1=2$.  By the Cartan-K\"ahler Theorem (see \cite{CFB}, Chapter 7)
 we conclude that
6-dimensional integral submanifolds, satisfying the independence condition,
exist through every point of $\scrU$, and that the construction of such
submanifolds depends on a choice of $2$ functions of one variable.

The additional 1-forms that generate the prolongation include
\begin{align*}
\psi_7 &= d(f_u) +\left(f_u^2+2\dfrac{f_u}{x+y}\right)dx+\left(f_u g_u-\dfrac{2}{(x+y)^2}\right)dy,\\
\psi_8 &= d(g_u) + \left(g_u^2 + 2\dfrac{g_u}{x+y}\right) dy +\left(f_u g_u - \dfrac{2}{(x+y)^2}\right)dx,
\end{align*}
which are actually defined on the original manifold $\R^{16}$.
These forms vanish on all integral elements of the original system, and if they had
been included in the ideal, it would have been involutive with $s_1 = 2$.

The vanishing of $\psi_7,\psi_8$ implies that $f_u$ and $g_u$ are
functions of $x$ and $y$ only.
Moreover, forms $\psi_7,\psi_8$ define a smaller Pfaffian system, involving only $f_u,g_u$ as functions of $x$ and $y$, which satisfies the Frobenius condition.  This means that $f_u(x,y)$ and $g_u(x,y)$ can be determined by solving systems of ODE.  Once these are determined, substituting the solutions into \eqref{b5fpx} and \eqref{b5gqy} gives a Frobenius system which may be solved for the functions $f_P(x)$ and $g_Q(y)$.  Then $f$ and $g$ may be determined by integrating first-order PDE, with $f$ including an arbitrary function of $x$ and $g$ an arbitrary function of $y$.

For example, by observing that $f_u + g_u$ must satisfy a Riccati equation as a function of $x+y$, we are led to a solution $$f_u = \dfrac{y}{(x+y)x}, \qquad g_u = \dfrac{x}{(x+y)y}.$$ Substituting these into \eqref{b5fpx},\eqref{b5gqy} leads to $f_P + g_Q = k(x+y)/(xy)$ for a constant $k$.  It is simplest to choose $k=0$ with $f_P = 1$ and $g_Q = -1$.  Integrating then gives the solution $$f =  P+\dfrac{y(u+Z)}{x(x+y)},\qquad g = -Q +\dfrac{x(u-Z)}{y(x+y)}.$$

\begin{prop} All \BT/s between (XIII) and the wave equation are holonomic.\end{prop}
\begin{proof}  As noted in \S\ref{zsec}, the holonomic condition is equivalent to
the Pfaffian system on $\B^6$ spanned by $dx,dy$ and
$$dp - f_u du \equiv (-g_Q f_u) dZ + f_P dP, \quad dq -g_u du\equiv (-f_P g_u) dZ + g_Q dQ, \mod dx, dy$$
being Frobenius.
It is straightforward to check that $d(f_P dP - g_Q f_u dZ)$ and $d(g_Q dQ - f_P g_u dZ)$ are
zero modulo $dx,dy$ and the 1-forms of the above EDS.  (For example,
$d(f_P) \equiv 0$ modulo $dx, dy, \psi_1, \ldots, \psi_8$, and
the same is true for $d(g_Q)$ and $d(f_u)$.)
\end{proof}

\bigskip
\subsubsection*{Equation IV}  This PDE has the form
$$s =\Fsigma(u,p,q) := \pm \dfrac{\alpha(p)\beta(q)}{u},$$
where $\alpha$ and $\beta$ are arbitrary solutions of the ODE $df/dt \pm t/f = K$ for some fixed $K\ne 0$.
(We will take the plus sign in these equations, the computation
for the other sign being completely analogous.)  In this case,
\eqref{genfzgz} takes the form
\begin{align}
f_Z = \left( \dfrac{\alpha(f)}u \left(K -\dfrac{g}{\beta(g)}\right) - f_u\right) g_Q, \label{G4fz} \\
g_Z = \left( \dfrac{\beta(g)}u \left(K -\dfrac{f}{\alpha(f)}\right) - g_u\right)f_P. \label{G4gz}
\end{align}
(From now on, instead of writing $\alpha(f)$ and $\beta(g)$,
$\alpha$ and $\beta$ will be understood to be composed with
$f(x,y,u,Z,P)$ and $g(x,y,u,Z,Q)$ respectively.)

Unlike in previous examples,
here it is not valid to conclude that $f$ and $g$ are linear in $P$ and $Q$.
In fact, differentiating
\eqref{G4fz} with respect to $Q$ gives
\begin{equation}\label{G4gqq}
0=(u \beta f_u -(K\beta-g)\alpha)g_{QQ}
+\dfrac{(\beta^2 -(K\beta - g)g)}{\beta^2}\alpha g_Q^2,
\end{equation}
enabling us to determine $g_{QQ}$.  (If the coefficient in front were identically
zero, then $\beta(q)$ would be identically equal to a constant times $q$,
which contradicts $K\ne 0$.)  Similarly, differentiating \eqref{G4gz} yields
\begin{equation}\label{G4fpp}
0 =(u \alpha g_u - (K\alpha -f)\beta) f_{PP}
+ \dfrac{(\alpha^2 -(K\alpha - f)f)}\beta f_P^2.
\end{equation}

The characteristic invariants for (IV) are
$$J_1 = \dfrac{r}{\alpha} - \dfrac{\alpha}{u}, \qquad J_2 = \dfrac{t}{\beta} - \dfrac{\beta}{u}.$$
Substituting for $p$ and $r$ from \eqref{pandq} and \eqref{randt} gives
\begin{equation}\label{B4J1after}
J_1 = \dfrac{f_x + f_u f + f_Z P + f_P R}{\alpha} - \dfrac{\alpha}{u},\
\end{equation}
so that $f_P /\alpha$ must be a function of $x$ and $P$ only.
Setting the derivatives of this with respect to $u,y$, and $Z$ equal to zero yields
\begin{equation}\label{G4fpu}
f_{Pu} = \left(K -\dfrac{f}{\alpha}\right)  \dfrac{f_P f_u}{\alpha}, \qquad
f_{Py} = \left(K -\dfrac{f}{\alpha}\right)  \dfrac{f_P f_y}{\alpha},\qquad
f_{PZ} = \left(K -\dfrac{f}{\alpha}\right)  \dfrac{f_P f_Z}{\alpha},
\end{equation}
where $f_y$ is given by \eqref{genfy} and $f_Z$ is given by \eqref{G4fz}.
Similarly, from the $T$ coefficient in $J_2$ we get that $g_Q/\beta$ must be a function of
$y$ and $Q$ only, and hence
\begin{equation}\label{G4gqu}
g_{Qu} = \left(K -\dfrac{g}{\beta}\right)  \dfrac{g_Q g_u}{\beta}, \qquad
g_{Qx} = \left(K -\dfrac{g}{\beta}\right)  \dfrac{g_Q g_x}{\beta},\qquad
g_{QZ} = \left(K -\dfrac{g}{\beta}\right)  \dfrac{g_Q g_Z}{\beta}.
\end{equation}
We may also differentiate \eqref{G4fz} and  \eqref{G4gz} to obtain
equations for $f_{Px}$ and $g_{Qy}$.

We encode the various first- and second-order partial differential
equations for $f$ and $g$ derived so far into an exterior differential system generated by 1-forms $\psi_1, \ldots, \psi_6$, as we did for equation (XIII).
Unlike the previous example,
we do not need to prolong, but instead obtain integrability conditions which take the form
\begin{equation}\label{G4intcond2}
\alpha g_u = \beta f_u,
\end{equation}
and
\begin{equation}\label{G4intcond3}
 (u \beta f_u -(K\beta-g)\alpha) g_Q = (u \alpha g_u - (K\alpha -f)\beta) f_P.
\end{equation}
(Note that, by using \eqref{G4fz}, \eqref{G4gz}, this implies that
$\alpha g_Z = \beta f_Z.$)
With these conditions incorporated into the EDS, it becomes involutive
with last nonzero character $s_1=2$.

Solutions of this system may be obtained by observing that the quantities
\begin{equation}\label{defgammadelta}
\lambda = \dfrac{f_u}{\alpha}, \quad \mu = \dfrac{f_Z}{\alpha},
\quad \gamma = \alpha - (\lambda f + \mu P)u,
\quad \delta = \beta - (\lambda g + \mu Q)u
\end{equation}
must be functions of $x,y,u$ and $Z$ only, and satisfy the following compatible
system of first-order PDE:
\begin{align*}
\dfrac{\di\lambda}{\di u} &=-u\lambda^3 + K \lambda^2 - \dfrac{2}{u}\lambda,
&
\dfrac{\di \mu}{\di u} &=-\dfrac{(u^2\lambda^2 - K u\lambda+1)}{u}\mu,
&
\dfrac{1}{\gamma}\dfrac{\di \gamma}{\di u} &=
\dfrac{1}{\delta}\dfrac{\di \delta}{\di u} = (K-u\lambda)\lambda,
\\
\dfrac{\di\lambda}{\di Z} &= \dfrac{\di\mu}{\di u},
&
\dfrac{\di\mu}{\di Z} &= (K-u\lambda) \mu^2,
&
\dfrac{1}{\gamma}\dfrac{\di \gamma}{\di Z} &=
\dfrac{1}{\delta}\dfrac{\di \delta}{\di Z} = (K-u\lambda)\mu,
\\
\dfrac{\di\lambda}{\di x} &= -\dfrac{(u^2\lambda^2 - K u\lambda+1)}{u^2}\gamma,
&
\dfrac{\di\mu}{\di x} &= (K-u\lambda) \dfrac{\gamma \mu}{u},
&
\dfrac{\di\delta}{\di x} &= (K - u\lambda) \dfrac{\gamma\delta}{u},
\\
\dfrac{\di\lambda}{\di y} &= -\dfrac{(u^2\lambda^2 - K u\lambda+1)}{u^2}\delta,
&
\dfrac{\di\mu}{\di y} &= (K-u\lambda) \dfrac{\delta \mu}{u},
&
\dfrac{\di \gamma}{\di y} &= \dfrac{\di\delta}{\di x}.
\end{align*}
A solution $(\gamma, \delta, \lambda,\mu)$ to this PDE system
may be constructed by integrating successively in the $u$-direction,
the $Z$-direction, the $x$-direction and the $y$-direction.  (Note that
the $x$-dependence of $\gamma$ and the $y$-dependence of $\delta$ are
given by arbitrary functions.)  Once $\gamma$ and $\delta$ are
known, they implicitly determine $f$ and $g$.

Using \eqref{holonomicxi}, one can check that
the resulting \BT/s are holonomic if and only if, in the above system,
$u^2\lambda^2 -Ku\lambda+1=0.$  Thus, holonomic transformations
exist, and depend on fewer arbitrary constants but the same number
of arbitrary functions.  For example, if $K=2$, then a solution to the above
system is given by
$$\lambda =\dfrac1{u}, \qquad \mu = \dfrac{-1}{Z -v(x)-w(y)},\qquad
\gamma = \dfrac{v'(x)u}{Z -v(x)-w(y)}, \qquad \delta = \dfrac{w'(y) u}{Z-v(x)-w(y)}.
$$
Then, using \eqref{defgammadelta},  a holonomic \BT/  is implicitly
defined by
$$\alpha(p) - p = \dfrac{(v'(x)-P)u}{Z-v(x)-w(y)},
\qquad
\beta(q) - q = \dfrac{(w'(y)-Q)u}{Z-v(x)-w(y)}.
$$

\subsection{Summary}
Besides equations (IV) and (XIII) discussed above,
we have also investigated the exterior differential
system for \BT/s to the wave equation for
equations (V), (VII), (IX), (XI) and (XII).
Even if explicit formulas are not available,
in each case we use the Cartan-K\"ahler Theorem to determine
(in terms of the last nonzero Cartan character) the size of the
solution set, in both the holonomic and non-holonomic cases.
The results are summarized in the table below.

\bigskip
\begin{tabular}{|c|c|c|c|}\hline
Equation & Monge-Integrable & Holonomic BTs & Non-holonomic BTs \\ \hline
I 	& no & yes, $s_1=2$	& yes, $s_1=2$ \\
II 	& no & no          &	 yes, $s_1=2$ \\
III 	& no & no 			& yes, $s_1=2$ \\
IV	& no & yes, $s_1=2$ & yes, $s_1=2$ \\
V	& no & yes, $s_1=2$ & yes, $s_1=2$ \\
VI	& yes& yes, $s_1=3$ & yes, $s_3=1$ \\
VII	& no & no & yes, $s_1=2$ \\
VIII	& yes& yes, $s_1=3$ & yes, $s_3=1$ \\
IX	& no & no			& yes, $s_1=2$ \\
X 	& yes& yes, $s_1=3$ & yes, $s_3=1$ \\
XI	& no  & yes, $s_1=2$ & no \\
XII	& no  & yes, $s_1=2$ & no \\
XIII	& no 	 & yes, $s_1=2$ & no \\ \hline
\end{tabular}

\medskip
Note that the approach described in \S\ref{solving}
is not feasible for the Monge-integrable equations (VI, VIII and X),
but for completeness we include them in the table,
together with the results from
the analysis in \S\ref{DtoBmonge}.  The Cartan character for
the system for holonomic \BT/s for such equations is variable,
depending on whether one considers the cases (1), (2), or (3), as
 described in \S\ref{DtoBmonge}.

It is interesting to note that equations (I), (IV) and (V) have both holonomic
and non-holonomic transformations, in roughly the same degree of generality.
In fact, it is possible that these two kinds of \BT/s linking the same pair
of equations may be closely related.  In our previous paper \cite{CI}, we pointed out that
the transformation (Z.I) is a composition of two simpler transformations, a holonomic
\BT/ to the wave equation, and a contact transformation from the wave equation to itself.
It is possible that, more generally, the non-holonomic transformations for these
equations are obtainable from holonomic transformations in this way.

\setcounter{section}{5}
\section{Concluding Remarks}\label{wconclusion}
In this section, we will indicate some interesting directions in which the results in
this paper might be extended, and some important questions about \BT/s to which the
techniques in this paper may be relevant.

\begin{enumerate}
\item  The set of equivalence classes (under contact transformations) of second-order \MA/ equations
to which the results of \S3 in this paper apply is relatively small, confined to the equations on
the Goursat-Vessiot list.  It would be interesting to see if the arguments in that section
could be applied to hyperbolic systems of class $k>1$.  In other words, given a hyperbolic system $\I$
of class $k$, linked to the standard wave equation by a \BT/, can one prove that the prolongation of $\I$ is
Darboux-integrable?
Likewise, given a hyperbolic EDS $\I$ of
class $k$, such that its prolongation is Darboux-integrable, does there exist a \BT/  between
$\I$ and the \MA/ system for the standard wave equation?  (The argument given at the end of
Chapter 7 in \cite{CFB} shows that there is a \BT/ between the wave equation and the {\em prolongation}
of $\I$; however, for practical purposes it is desirable to have
a \BT/ between systems of as low an order as possible, so that one has a smaller
system of ODE to solve in order to construct solutions.)  These hyperbolic systems would include, for example, the \MA/ equations which
are Darboux-integrable at third order, which have not been classified and are thought
to comprise a much larger set.

\item
It is a theorem of Sophus Lie that no \MA/ equation of the form $u_{xy}=f(u)$ is Darboux-integrable
(after arbitrary many prolongations) except when $f(u) = \exp(a u +b)$ for constants $a$ and $b$ (see \cite{GLecons}, Chapter IX).
Consequently, important equations like sine-Gordon cannot have a \BT/ to the wave equation.
Instead, the \BT/ \eqref{SGsingleBT} for sine-Gordon produces solutions to the {\em same} PDE
as we started with.  This is known as an {\em auto-\BT/};\footnote{This terminology is not universally accepted; Hongyou Wu \cite{Wu} has proposed
that transformations between different PDEs be known as {\em Miura transformations},
and the term \BT/ be reserved for what we are calling auto-\BT/s.}
 such transformations play an important
role in the theory of completely integrable PDE \cite{Fordy}.

It is therefore of interest to
try to identify those \MA/ equations which have non-trivial auto-\BT/s.  We remark that
for such transformations, the \MA/ systems $\I$ on $\scM$ and $\Ibar$ on $\Mbar$
must be contact-equivalent, i.e.,. there must be a diffeomorphism $\Phi: \scM \to \Mbar$ which
pulls back $\Ibar$ to $\I$.

\setlength{\unitlength}{2pt}
\begin{center}
\begin{picture}(40,30)(0,0)
\put(5,5){\makebox(0,0){$\scM$}}
\put(35,5){\makebox(0,0){$\Mbar$}}
\put(20,25){\makebox(0,0){$\B$}}
\put(16,21){\vector(-3,-4){9}}
\put(24,21){\vector(3,-4){9}}
\put(10,5){\vector(1,0){21}}
\put(18,7){$\Phi$}
\put(7,16){$\pi$}
\put(31,16){$\pibar$}
\end{picture}
\end{center}
Necessary conditions for the existence of such a diffeomorphism may be derived
from the fact that it is required to preserve the differential invariants of the \MA/ systems.
(See \cite{BGG}, \S2.1, for a derivation of these invariants using the method of equivalence.)

\item Our previous paper \cite{CI} began the exploration of {\em parametric} \BT/s using
the method of equivalence.   Such transformations contain an arbitrary parameter in the
\BB/ system; for example, an arbitrary nonzero parameter $\lambda$ may be interpolated
in the sine-Gordon auto-\BT/  \eqref{SGsingleBT} to give
$$
\begin{aligned}
v_x - u_x&=  \dfrac{\lambda}{2}\sin((u+v)/2), \\
v_y + u_y&=  -\dfrac{1}{2\lambda}\sin((u-v)/2).
\end{aligned}
$$
One observes that this system differs from \eqref{SGsingleBT} merely by scaling
$x$ by $\lambda$ and $y$ by $\lambda^{-1}$---a change of variables which is a
symmetry of the sine-Gordon equation but not of the system \eqref{SGsingleBT}.  This scaling symmetry
can also be applied to the \BT/ \eqref{Liouvillewave}, to produce a parametric transformation
$$
z_x = u_x - 2\lambda\exp((u+z)/2), \qquad
z_y = -u_y +\dfrac{1}{\lambda}\exp((u-z)/2),
$$
where $u(x,y)$ satisfies Liouville's equation and $z(x,y)$ solves the wave equation.
(In fact, this transformation is derived from the most general form \eqref{Liouvillegeneral}
by setting $k=1$ and choosing $v(x)=2\ln \lambda$ and $w(y)=0$.)  In \cite{CI} it is shown that
these transformations can be generated from a non-parametric \BT/  by starting with a symmetry vector field
on $\scM$, choosing a lift into $\B^6$ which is {\em not} a symmetry of the Pfaffian system
$\J$, but such that pulling $\J$ back to $\B \times \R$ via the 1-parameter
family of diffeomorphisms generated by the lift gives a family of transformations.
The same approach can be taken with other transformations discussed in \S4.1; for example,
the transformation \eqref{G2bt} may be generalized to a parametric transformation
$$z z_x - \lambda u u_x = \sqrt{z^2 - \lambda u^2}\sqrt{\lambda+z_x^2}, \qquad
z z_y - \lambda u u_y = \sqrt{z^2 - \lambda u^2}\sqrt{\lambda+z_y^2}
$$
for $\lambda >0$, where $u$ satisfies (II) and $z$ solves the wave equation.
(This is obtained by starting with the symmetry of (II) that simultaneously
scales $u$, $x$ and $y$.)

With these examples in evidence, and given the importance of parametric \BT/s
in the study of `soliton' equations, it is desirable to characterize those transformations
that may be made to depend on an arbitrary parameter by lifting symmetry vector fields.

\end{enumerate}

\end{document}